\documentclass[reqno,twoside,11pt]{amsart}
 
\hoffset - 1.5cm

\usepackage{graphicx}

\usepackage{amstext}

\usepackage[T1]{fontenc}

\usepackage{color}
\usepackage{amsmath}
\usepackage{amsthm}
\usepackage{amssymb}
\usepackage[mathscr]{euscript}
\usepackage[colorinlistoftodos,prependcaption,textsize=tiny]{todonotes}
\usepackage{enumerate}

\newtheorem{Theorem}{Theorem}[section]
\newtheorem{Fact}[Theorem]{Fact}
\newtheorem{Lemma}[Theorem]{Lemma}
\newtheorem{Proposition}[Theorem]{Proposition}
\newtheorem{Corollary}[Theorem]{Corollary}
\theoremstyle{definition}
\newtheorem{Definition}[Theorem]{Definition}
\newtheorem{Example}[Theorem]{Example}
\newtheorem{Remark}[Theorem]{Remark}

\newcommand{\numsec}{
\setcounter{Theorem}{0}
\renewcommand\theequation{\arabic{section}.\arabic{equation}}
\renewcommand\theTheorem{\arabic{section}.\arabic{Theorem}}
\renewcommand\theDefinition{\arabic{section}.\arabic{Theorem}}
\renewcommand\theRemark{\arabic{section}.\arabic{Theorem}}
\renewcommand\theLemma{\arabic{section}.\arabic{Theorem}}
\renewcommand\theFact{\arabic{section}.\arabic{Theorem}}
\renewcommand\theProposition{\arabic{section}.\arabic{Theorem}}
\renewcommand\theCorollary{\arabic{section}.\arabic{Theorem}}
\renewcommand\theExample{\arabic{section}.\arabic{Theorem}}
\renewcommand\theProperty{\arabic{section}.\arabic{Theorem}}
}

\newcommand{\numsubsec}{
\setcounter{Theorem}{0}
\renewcommand\theequation{\arabic{section}.\arabic{subsection}.\arabic{equation}}
\renewcommand\theTheorem{\arabic{section}.\arabic{subsection}.\arabic{Theorem}}
\renewcommand\theDefinition{\arabic{section}.\arabic{subsection}.\arabic{Theorem}}
\renewcommand\theRemark{\arabic{section}.\arabic{subsection}.\arabic{Theorem}}
\renewcommand\theLemma{\arabic{section}.\arabic{subsection}.\arabic{Theorem}}
\renewcommand\theFact{\arabic{section}.\arabic{subsection}.\arabic{Theorem}}
\renewcommand\theProposition{\arabic{section}.\arabic{subsection}.\arabic{Theorem}}
\renewcommand\theCorollary{\arabic{section}.\arabic{subsection}.\arabic{Theorem}}
\renewcommand\theExample{\arabic{section}.\arabic{subsection}.\arabic{Theorem}}
\renewcommand\theProperty{\arabic{section}.\arabic{subsection}.\arabic{Theorem}}
}

\newcommand{\ba}{\begin{array}}
\newcommand{\bc}{\begin{center}}
\newcommand{\bd}{\begin{description}}
\newcommand{\bdm}{\begin{displaymath}}
\newcommand{\be}{\begin{enumerate}}
\newcommand{\beq}{\begin{equation}}
\newcommand{\bdf}{\begin{Definition}}
\newcommand{\bex}{\begin{Example}}
\newcommand{\bft}{\begin{Fact}}
\newcommand{\bl}{\begin{Lemma}}
\newcommand{\bp}{\begin{Proposition}}
\newcommand{\br}{\begin{Remark}}
\newcommand{\bt}{\begin{Theorem}}
\newcommand{\bco}{\begin{Corollary}}
\newcommand{\bh}{\begin{Hipothesis}}
\newcommand{\ea}{\end{array}}
\newcommand{\ec}{\end{center}}
\newcommand{\ed}{\end{description}}
\newcommand{\edm}{\end{displaymath}}
\newcommand{\ee}{\end{enumerate}}
\newcommand{\eeq}{\end{equation}}
\newcommand{\edf}{\end{Definition}}
\newcommand{\eex}{\end{Example}}
\newcommand{\eft}{\end{Fact}}
\newcommand{\el}{\end{Lemma}}
\newcommand{\ep}{\end{Proposition}}
\newcommand{\er}{\end{Remark}}
\newcommand{\et}{\end{Theorem}}
\newcommand{\eco}{\end{Corollary}}
\newcommand{\eh}{\end{Hipothesis}}

\newcommand{\bC}{\mathbb{C}}

\newcommand{\bE}{\mathbb{E}}

\newcommand{\bH}{\mathbb{H}}
\newcommand{\bI}{\mathbb{I}}

\newcommand{\bN}{\mathbb{N}}

\newcommand{\bP}{\mathbb{P}}

\newcommand{\bR}{\mathbb{R}}

\newcommand{\bV}{\mathbb{V}}
\newcommand{\bW}{\mathbb{W}}
\newcommand{\bX}{\mathbb{X}}
\newcommand{\bY}{\mathbb{Y}}

\newcommand{\cC}{\mathcal{C}}

\newcommand{\cE}{\mathcal{E}}
\newcommand{\cF}{\mathcal{F}}

\newcommand{\cH}{\mathcal{H}}
\newcommand{\cI}{\mathcal{I}}

\newcommand{\cN}{\mathcal{N}}
\newcommand{\cO}{\mathcal{O}}

\newcommand{\cR}{\mathcal{R}}

\newcommand{\cT}{\mathcal{T}}

\newcommand{\im}{\mathrm{ im \;}}
\newcommand{\diag}{\mathrm{ diag \;}}
\numberwithin{equation}{section} \errorcontextlines=0

\newcommand{\vp}{\varphi}

\newcommand{\morse}{\mathrm{m^-}}

\newcommand{\sone}{S^1}
\newcommand{\sotwo}{SO(2)}

\newcommand{\ds}{\displaystyle}

\newcommand{\sptwon}{Sp(2N,\bR)}

\newcommand{\subg}{\overline{\mathrm{sub}}(G)}
\newcommand{\subgc}{\overline{\mathrm{sub}}[G]}
\newcommand{\subh}{\overline{\mathrm{sub}}(H)}

\newcommand{\gcw}{G\text{-CW}}
\newcommand{\gls}{$G$-$\mathcal{LS}$}

\newcommand{\chig}{\chi_G}

\newcommand{\upsg}{\Upsilon_G}
\newcommand{\upsh}{\Upsilon_H}
\newcommand{\upssone}{\Upsilon_{\sone}}
\newcommand{\ci}{\cC\cI}
\newcommand{\cig}{\ci_G}

\newcommand{\scig}{\mathscr{CI}_G}
\newcommand{\sci}{\mathscr{CI}}

\newcommand{\lpm}{\lambda_{\pm}}

\newcommand{\brtwon}{\bR^{2N}}

\begin{document}

\title[Periodic solutions of symmetric Hamiltonian systems]{Periodic solutions of symmetric Hamiltonian systems}

\author{Daniel Strzelecki}
\address{Faculty of Mathematics and Computer Science\\
Nicolaus Copernicus University \\
PL-87-100 Toru\'{n} \\ ul. Chopina $12 \slash 18$ \\
Poland}

\email{danio@mat.umk.pl}

\date{\today}


\begin{abstract} 
This paper is devoted to the study of periodic solutions of Hamiltonian system $\dot z(t)=J \nabla H(z(t))$, where $H$ is symmetric under an action of a compact Lie group. We are looking for periodic solutions in a nearby of non-isolated critical points of $H$ which form orbits of the group action. We prove Lyapunov-type theorem for symmetric Hamiltonian systems.
\end{abstract}

\maketitle

\noindent\textbf{Note on this version.} This is a corrected version of the paper published in Arch. Ration. Mech. Anal. 237 (2020), 921--950. It incorporates the corrections from the author's Corrigendum: the reasoning leading to assumption \textbf{(A7.5)} in Section~5 has been corrected, assumption \textbf{(A7.5)} has been replaced by a stronger one, and Theorem~\ref{T:intro} in the Introduction, Theorem~\ref{T:cor3} and its application in Section~6 have been adjusted accordingly. The main result of the paper (Theorem~\ref{T:general}) and the methodology are not affected.

\section{Introduction}

Consider a first-order system 
\beq\label{hamsys}\tag{HS}
\dot z(t)=J \nabla H(z(t))
\eeq 
on $\brtwon$, where $J=\left[\ba{rr} 0 & I \\ -I & 0 \ea \right]$ is the standard symplectic matrix and $H:\brtwon\to\bR$ is a Hamiltonian of the class $C^2$.

The existence of periodic orbits in Hamiltonian dynamics is an important and widely studied problem. In 1895 Lyapunov \cite{LIAPUNOV} proved his center theorem i.e. the existence of one-parameter family of periodic solutions of \eqref{hamsys} tending to a non-degenerate equilibrium. The next important result of Weinstein \cite{WEINSTEIN} shows the existence of at least $N$ geometrically distinct periodic solutions at any energy level of the Hamiltonian $H$. The further development of Weinstein theorem was performed by Moser \cite{MOSER}. In 1978 Fadell and Rabinowitz proved the lower bound for the number of small nontrivial solutions of \eqref{hamsys} depending on the period. See \cite{RAB6} for the general overview of the results up to 1982. The results of Weinstein and Moser were generalized by Bartsch in 1997, \cite{BARTSCH2}. The problem of the existence of periodic solutions of \eqref{hamsys} in a case of degenerate equilibrium was also studied by Szulkin \cite{SZULKIN} and Dancer with Rybicki \cite{DARY} who generalized classical result of Lyapunov.

Suppose now that the compact Lie group $\Gamma$ acts unitary on $\brtwon$ and $H \in C^2(\bR^{2N},\bR)$ is a $\Gamma$-invariant potential i.e. $H(\gamma z)=H(z)$ for any $\gamma\in\Gamma$ and $z\in\brtwon$. The study of the existence of periodic solutions in this case was performed by Montaldi, Roberts and Stewart \cite{MOROST} who had proved equivariant version of Weinstein-Moser theorem. In 1993  Bartsch \cite{BARTSCH} generalized the theorem of Montaldi, Roberts and Stewart for the wider class of a group actions which allows him to generalize the result of Fadell and Rabinowitz also. However, the authors mentioned above have assumed that critical point $z_0$ of $H$ is a fixed point of group action i.e. the orbit of action consists of one point. Then $z_0$ can be an isolated critical point. 

We study a more general case. Assume that $z_0$ is a critical point of Hamiltonian $H$. Since $H$ is $\Gamma$-invariant, $\Gamma(z_0)=\{\gamma z_0\,:\,\gamma\in\Gamma\}\subset (\nabla H)^{-1}(0)$, i.e. the orbit $\Gamma(z_0)$ consists of critical points of $H$ and, therefore, stationary solutions of the equation \eqref{hamsys}, see Remark \ref{eqgrad}. Hence, if $\dim \Gamma_{z_0}<\dim \Gamma$ then the orbit is at least an one-dimensional manifold and, as a consequence, critical points are not isolated. Therefore the results mentioned above are not applicable. 
We are going to prove sufficient conditions for the existence of non-constant periodic solutions of an autonomous Hamiltonian system in the presence of symmetries of a compact Lie group, the problem \eqref{hamsys}, in any neighborhood of the orbit $\Gamma(z_0)$, see the main result Theorem \ref{T:general} and Theorems \ref{T:cor1}, \ref{T:cor2}, \ref{T:cor3}.

This article is organized as follows. In section \ref{preliminaria} we recall some basic definitions and notions of group theory and equivariant topology. Equivariant Conley index which is a main tool of our reasoning is shortly defined in subsection \ref{conleyIndexSection}. Furhter, we recall the notion of Euler ring, equivariant Euler characteristic and its generalization, see Remarks \ref{propOfChi}, \ref{stabilization of upsilon}. In Theorem \ref{smash theorem for sci} we recall the very important theorem connecting equivariant Euler characteristic, equivariant Conley index and the idea of orthogonal section introduced in the paper \cite{PCHRYST}.  In subsection 2.5 we formulate so called equivariant splitting lemma - the theorem which allows us to simplify the study of Conley indexes up to the linear case in Lemma \ref{main-lemma0}.

In section \ref{variationa} we parameterize the equation \eqref{hamsys} to study the solutions with constant period $2\pi$ in the equation \eqref{hamsys2}. Next we introduce appropriate Sobolev space $\bE$, the action of the group $\Gamma\times\sone$ on it and we define variational functional $\Phi:\bE\to\bR$ (see formulas \eqref{PhiDef0}, \eqref{PhiDef}) such that $2\pi$-periodic solutions of the system \eqref{hamsys2} are in bijective correspondence with $\sone$-orbits of critical points of $\Phi$. In this way we begin to study the equation \eqref{gradeq}.
Further, we analyze the linear Hamiltonian system \eqref{hamsyslin}; it is a base for the last step in the proof of the main result of the paper.

Section \ref{S:main} is devoted to the formulation and the proof of the main result of this paper, Theorem \ref{T:general}. The notion of bifurcation theory is recalled in Definition \ref{D:globBif} and in the nearby text. In Theorems \ref{necCond}, \ref{sufficient} we formulate the necessary and sufficient condition for the existence of global bifurcation of solutions of the equation \eqref{gradeq}. The last part of this section is devoted to the proof of the change of equivariant Euler characteristics of equivariant Conley indexes i.e. the formula \eqref{suffCondition}. 
Firstly, we reduce our task to the space orthogonal to the orbit, see Lemma \ref{sectionConditionLemma} and the text above them. Next in Lemma \ref{main-lemma0} we reduce the problem to the linear case. To study it we prove Theorem \ref{lastLemma}. To finish the proof of the main result we study the minimal periods and convergence of new solutions in Remarks \ref{R:final1}, \ref{R:final2}.

In the fifth section we reformulate the main result to make the assumptions easier to verify. The most friendly version of our result is the following theorem (see Theorem \ref{T:cor3}).
\bt\label{T:intro}
Let $H:\brtwon\to\bR$ be a $\Gamma$-invariant Hamiltonian of the class $C^2$. Let $z_0$ be a critical point of $H$ such that $\Gamma_{z_0}=\{e\}$ and the orbit $\Gamma(z_0)$ is isolated in $(\nabla H)^{-1}(0)$. Assume that $m^+(\nabla^2 H(z_0))> N$ or $m^-(\nabla^2 H(z_0))> N$ and $\deg(\nabla H_{\mid T^{\perp}_{z_0} \Gamma(z_0)},B(z_0,\epsilon),0) \neq 0$ for sufficiently small $\epsilon$. Then there exists a connected family of non-stationary periodic solutions of the system $\dot z(t)=J \nabla H(z(t))$ emanating from the stationary solution $z_0$ such that periods (not necessarily minimal) of solutions in the small neighborhood of $z_0$ are close to $ 2\pi \slash \beta_{j}$, where $i\beta_j$, $\beta_j>0$, is some eigenvalue of $J\nabla^2 H(z_0)$.
\et 
\noindent For the two other versions see Theorems \ref{T:cor1} and \ref{T:cor2}.

Furhter, we show that Lyapunov-type theorem of Dancer and Rybicki (Theorem \ref{T:DARY}) is generalized by the main result of this paper - Theorem \ref{T:general}. In the last part of this section we reformulate the second-order Newtonian system to the Hamiltonian one. Then the two symmetric versions of the Lyapunov center theorem, Theorem \ref{SLCT} proven in \cite{PCHRYST} and Theorem \ref{SLCTM} proven in \cite{PCHRYST2}, are also the consequences of the results proven in this paper.

The last section is devoted to an interesting application of the abstract results of this paper. We study the existence of quasi-periodic motions of the satellite in a nearby of a geostationary orbit of an oblate spheroid. In order to do this we consider a gravitational motion in the rotating frame where the corresponding Hamiltonian is given by formula \eqref{E:hamApp}. It is $\sotwo$ invariant and possesses a critical point which represents the geostationary orbit in the original coordinates. Theorem \ref{T:cor3} will be directly applied in this problem to prove the existence of trajectories with arbitrarily small deviations from the geostationary ones.

\section{Preliminaries}\label{preliminaria}
In this section we recall the basic material on equivariant topology from \cite{DIECK, KBO} and prove some preliminary results. Throughout this section $G$ stands for a compact Lie group.

\subsection{Groups and their representations}
Denote  by $\subg$ the set of all closed subgroups of  $G$. We say that two subgroups $H, H' \in \subg$ are conjugate in $G$ if there exists $g \in G$ such that $H=gH'g^{-1}.$ The conjugacy is an equivalence relation on $\subg.$ The class of $H \in \subg$ will be denoted by $(H)_G$ and the set of conjugacy classes we denote by $\subgc$. 

If $x \in \bR^n$ then $G(x)=\{gx : g \in G\}$ is the orbit through $x$ and a group $G_x=\{g \in G : g x = x\} \in \subg$ is called the isotropy group of $x$. The isotropy groups of the elements of common orbit are conjugate i.e if $G(x_1)=G(x_2)$ then $(G_{x_1})_G=(G_{x_2})_G$.  An open subset $\Omega \subset \bR^n$ is said to be $G$-invariant if  $G(x) \subset \Omega$ for every $x \in \Omega$. Note that the orbit $G(x)$ is a smooth $G$-manifold which is $G$-diffeomorphic to $G \slash G_x$.

Below we recall the notion of an admissible pair, which was introduced in \cite{PCHRYST}, where one can find some examples and properties.

\bdf \label{admissible}
Fix $H \in \subg.$ A pair $(G,H)$ is said to be \textit{admissible} if for any $K_1,K_2 \in \subh$  the following condition is satisfied: 
$\text{ if } (K_1)_H \neq (K_2)_H \text{ then } (K_1)_G \neq (K_2)_G.$
\edf

Note that if $\Gamma$ is a compact Lie group, then the pair $(\Gamma \times \sone,\{e\} \times \sone)$ is admissible, see Lemma 2.1 of \cite{PCHRYST}. This property will play a crucial role in the proof of the main result, Theorem \ref{T:general}.

Recall that a unitary group $U(N)$ id defined by
\[
U(N)=Sp(2N,\bR)\cap O(2N)
\]
where 
\[
\sptwon=\{A\in M_{2N\times 2N}(\bR)\,:\, A^T J A=J\}
\]
is a symplectic group and
\[
O(2N,\bR)=\{A\in M_{2N\times 2N}(\bR)\,:\, A^TA=Id\}
\]
is an orthogonal group. In particular if $A\in U(N)$ then $JA=AJ$. 
Note that $U(N)$ is a compact subgroup of $GL(2N)$, $A\in U(N)$ implies $A^T=A^{-1}\in U(N)$ and $|\det A|=1$.

Let $\rho : G \to U(N)$  be a continuous homomorphism. The space $\brtwon$
with the $G$-action defined by $G \times \brtwon \ni (g,x) \to \rho(g)x \in \brtwon$ we call a real,  unitary representation of $G$.  To simplify notation we write $gx$ instead of $\rho(g)x$ and $\brtwon$ instead of $\bV$ if the homomorphism is known.
 
Two unitary representations of $G$, say $\bV=(\brtwon,\rho), \bV'=(\brtwon,\rho'),$ are equivalent (briefly $\bV \approx_G \bV'$) if there exists an equivariant linear isomorphism $L : \bV \to \bV'$ i.e. the isomorphism $L$ satisfying $L(gx)=gL(x)$ for any $g \in G, x \in \brtwon.$ 
Put $D(\bV)=\{x \in \bV : \| x \| \leq 1\},\, S(\bV)=\partial D(\bV)$, $S^{\bV}=D(\bV) \slash S(\bV)$ and $B_{r}(\bV)=\{x \in \bV : \| x \| < r\}$. Since the representation $\bV$ is orthogonal in particular, these sets are $G$ invariant.

\subsection{Equivariant maps}
Let  $(\bV,\langle\cdot,\cdot\rangle)$ be a unitary $G$-representation. Fix a $G$-invariant open  subset $\Omega\subset\bV$.

\bdf
We say $\phi : \Omega \to \bR$ of class $C^k$ is \textit{ $G$-invariant $C^k$-potential}, if $\phi(gx)=\phi(x)$ for every $g \in G$ and $x \in \Omega$. The set of $G$-invariant $C^k$-potentials will be denoted by $C^k_{G}(\Omega,\bR)$.
\edf

\bdf 
A map $\psi : \Omega \to \bV$ of the class $C^{k-1}$ is called \textit{$G$-equivariant $C^{k-1}$-map}, if $\psi(gx)=g \psi(x)$ for every $g \in G$ and  $x \in \Omega.$ Then we write $\psi\in C^{k-1}_{G}(\Omega,\bV).$
\edf

 Fix $\vp \in C^2_G(\Omega,\bR)$. By $\nabla \vp, \nabla^2 \vp$ we denote the gradient and the Hessian of $\vp$, respectively. $\morse(A)$ denotes the Morse index of symmetric matrix $A$  i.e. the sum of the  multiplicities of negative eigenvalues of $A$. Similarly by the $m^+(A)$ we denote the number of positive eigenvalues of $A$ counted with multiplicities.

\br \label{eqgrad} It is known that if $\vp \in C^k_G(\Omega,\bR)$ then $\nabla \vp \in C^{k-1}_G(\Omega,\bV)$. Since $\nabla\vp$ is $G$-equivariant if $x_0 \in (\nabla \vp)^{-1}(0)$ then $G(x_0)   \subset (\nabla \vp)^{-1}(0)$  i.e. critical points form orbits of a group action.  If $\nabla \vp(x_0)=0$ then  $\nabla \vp(\cdot)$ is fixed on $G(x_0)$ and therefore $T_{x_0} G(x_0) \subset\ker \nabla^2 \vp(x_0) $. As a consequence $\dim \ker \nabla^2 \vp(x_0) \geq \dim T_{x_0} G(x_0)=\dim G(x_0).$ 
\er

\subsection{Equivariant Conley index} \label{conleyIndexSection}
In this section we shortly recall the construction of equivariant Conley index introduced by Izydorek \cite{IZYDOREK}, see also \cite{GEIZPR,RYBAKOWSKI}.

Denote by $\cF_{\ast}(G)$ the category of finite pointed $\gcw$-complexes (see \cite{DIECK} for definition and examples). The $G$-homotopy type of $\bX\in\cF_{\ast}(G)$ we denote by $[\bX]_G\in \cF_{\ast}[G]$ (or $[\bX]$ when no confusion can arise) and by $\cF_{\ast}[G]$ the set of $G$-homotopy types of elements of $\cF_{\ast}(G)$. If $\bX$ is a $G$-CW-complex without a base point, then we denote by $\bX^+$ a pointed $G$-CW-complex $\bX^+=\bX\cup\{\ast\}$. By $\cig(S,\vartheta)$ we denote a finite dimensional $G$-equivariant Conley index of an isolated invariant set $S$ under a $G$-equivariant vector field $\vartheta$, see \cite{BARTSCH,FLOER,GEBA,SMWA} for the definition. Recall that $\cig(S,\vartheta)\in\cF_{\ast}[G]$.

Let $\xi=(\bV_n)_{n=0}^{\infty}$ be a sequence of finite-dimensional orthogonal $G$-representations.
\bdf\label{G-spectrum}
A pair $\cE(\xi)=\left((\cE_n)_{n=n(\cE(\xi))}^{\infty}, (\varepsilon_n)_{n=n(\cE(\xi))}^{\infty}\right)$, where $n(\cE(\xi))\in\bN$, is called a $G$-spectrum of type $\xi$ if
\begin{enumerate}
\item $\cE_n\in\cF_{\ast}(G)$ for $n\geq n(\cE(\xi))$,
\item $\varepsilon_n\in Mor_G(S^{\bV_n}\wedge \cE_n,\cE_{n+1})$ for $n\geq n(\cE(\xi))$,
\item there exists $n_1(\cE(\xi))\geq n(\cE(\xi))$ such that for $n>n_1(\cE(\xi))$, $\varepsilon_n$ is a $G$-homotopy equivalence.
\end{enumerate}
\edf

 The set of $G$-spectra of type $\xi$ is denoted by $GS(\xi)$. We can also define $G$-homotopy equivalence of two spectra $\cE(\xi),\cE'(\xi)$ (see \cite{IZYDOREK} for the details). The $G$-homotopy type of a $G$-spectrum $\cE(\xi)$ we denote by $[\cE(\xi)]_G$ (or shorter $[\cE(\xi)]$) and the set of $G$-homotopy types of $G$-spectra by $[GS(\xi)]$ or simply $[GS]$ when $\xi$ is fixed.

\br
It follows from definition that the $G$-homotopy type $[\cE(\xi)]$ of spectrum $\cE(\xi)=\left((\cE_n)_{n=n(\cE(\xi))}^{\infty}, (\varepsilon_n)_{n=n(\cE(\xi))}^{\infty}\right)$ depends only on the sequence $(\cE_n)_{n=n_1(\cE(\xi))}^{\infty}$.
\er

Let $(\bH,\langle\cdot,\cdot\rangle)$ be an infinite-dimensional orthogonal Hilbert representation of a compact Lie group $G$. Let $L:\bH\to\bH$ be a linear, bounded, self-adjoint and $G$-equivariant operator such that

\be \label{B-properties}
\item [(B.1)] $\bH=\overline{\bigoplus_{n=0}^{\infty} \bH_n}$, where all subspaces $\bH_n$ are mutually orthogonal G-representations of finite dimension,
\item [(B.2)]$H_0=\ker L$ and $L(\bH_n)=\bH_n$ for all $n\geq 1$,
\item [(B.3)]$0$ is not an accumulation point of $\sigma (L)$.
\ee

Put $\bH^n:=\bigoplus_{k=0}^n\bH_k$, denote by $P_n:\bH\to\bH^n$ the orthogonal projection onto $\bH^n$ and define the subspace of $\bH_k$ corresponding to the positive part of spectrum of $L$ by $\bH_k^+$.
Consider a functional $\Phi:\bH\to\bR$ such that $\nabla\Phi(x)=Lx+\nabla K(x)$, where $\nabla K\in C^1_H(\bH,\bH)$ is completely continuous. Denote by $\vartheta$ a \gls-flow, see Definition 2.1 of \cite{IZYDOREK}, generated by $\nabla\Phi$. Let be $\mathcal{O}$ an isolating $G$-neighborhood for $\vartheta$ and put $\cN=Inv_{\vartheta}\cO$. Set $\xi=(\bH_k^+)_{k=1}^{\infty}$. Let $\Phi_n:\bH^n\to\bR$ be given by $\Phi_n=\Phi_{\mid \bH^n}$ and $\vartheta_n$ denotes the $G$-flow generated by $\nabla\Phi_n$. Note that $\nabla\Phi_n(x)=Lx+P_n\circ\nabla K(x)$. Choose $n_0$ such that for $n\geq n_0$ the set $\cO_n:=\cO\cap\bH^n$ is an isolating $G$-neighborhood for the flow $\vartheta_n$. Then the set $Inv_{\vartheta_n}(\cO_n)$ admits a $G$-equivariantindex pair $(Y_n,Z_n)$. 

We define a spectrum $\cE(\xi):=(Y_n/Z_n)_{n= n_0}^{\infty}$. 
Then the equivariant Conley index of $\cO$ with respect to the flow $\vartheta$ is given by
$
\scig(\cO,\vartheta):=[\cE(\xi)]\in [GS].
$
Since isolated invariant set $\cO$ is defined by isolating neighborhood $\cN$ and the flow is related to vector field $\nabla\Phi$ we will also write $\scig(\cN,\nabla\Phi)$.

\subsection{Equivariant Euler characteristic}

Let $\left(U(G),+,\star)\right)$ be the Euler ring of  $G$, see 
\cite{DIECK} for the definition and more details. 
Let us briefly recall that the Euler ring $U(G)$ is commutative, generated by $\chig({G/H}^+)$, where $(H)\in\subgc$ with the unit $\bI_{U(G)}=\chig({G/G}^+)$, where $\chig:\cF_{\ast}[G]\to U(G)$ is the universal additive invariant for finite pointed $G$-CW-complexes known as the equivariant Euler characteristic.

\br\label{propOfChi}
Below we present some properties of the Euler characteristic $\chi_G(\cdot)$.
\begin{itemize}
\item  For  $\bX, \bY \in \cF_{\ast}(G)$ we have: $\chi_G(\bX)+\chi_G(\bY)=\chi_G(\bX\vee\bY)$ and $\chi_G(\bX)\star\chi_G(\bY)=\chi_G(\bX\wedge\bY)$.

\item If $\bW$ is a $G$-representation then $\chi_G(S^{\bW})$ is an invertible element of $U(G)$, see \cite{DGR}.

\item If $\bV,\bV'$ are $G$-representations such that $\dim \bV>\dim \bV'$ but $\bV\not\approx_G \bV'\oplus\bW$, where $\bW$ is even-dimensional trivial $G$-representations then\[
\chig(S^{\bV})\neq\chig(S^{\bV'}).
\]
For the prove of this fact see Lemma 3.4 in \cite{GLG}.
\end{itemize}
\er

 There is a natural extension of the equivariant Euler characteristic for finite pointed $G$-CW-complexes to the category of $G$-equivariant spectra due to Go{\l}\c{e}biewska and Rybicki \cite{GORY1}. 

 Let $\xi=(\bV_n)_{n=0}^{\infty}$ and put $\bV^n=\bV_0\oplus\bV_1\oplus\ldots\oplus\bV_n$, for $n\geq 0$. Recall that due to Remark \ref{propOfChi} an element $\chig(S^{\bV^n})$ is invertible in the Euler ring $U(G)$ and define a map $\upsg:[GS(\xi)]\to U(G)$ by the following formula
\beq \label{gchargspe}
\upsg([\cE(\xi)])=\lim_{n\to\infty} \left(\chig\left(S^{\bV^{n-1}}\right)^{-1}\star \chig(\cE_n)\right).
\eeq

\br\label{stabilization of upsilon}
It was shown in \cite{GORY1}  that $\upsg$ is well-defined. In fact 
\beq
\upsg([\cE(\xi)])=\chig\left(S^{\bV^{n_1(\cE)-1}}\right)^{-1}\star \chig(\cE_{n_1(\cE)})),
\eeq
where $n_1(\cE)=n_1(\cE(\xi))$ comes from Definition \ref{G-spectrum}.
\er

\br\label{Upsilon on finite complex}
Note that a finite pointed $G$-CW-complex $\bX$ can be considered as a \textit{constant} spectrum $\cE(\xi)$, where $\cE_n=\bX$ for all $n\geq 0$ and $\xi$ is a sequence of trivial, one-point representations. Then
$\ds
\upsg([\bX])=\upsg([\cE(\xi)])={\bI_{U(G)}}^{-1}\star\chig([\bX])=\chig([\bX]).
$
Therefore we can treat $\scig$ and $\upsg$ as natural extensions of $\cig$ and $\chig$ respectively.
\er

 By Theorems 3.1, 3.5  of \cite{GORY1} we obtain the following product formula.

\bt\label{product formula}
If $\cN_1,\,\cN_2$ are isolated $G$-invariant sets for the local $G$-$\mathcal{LS}$ flows generated by $\nabla\Psi_1$ and $\nabla\Psi_2$ respectively then
\[
\upsg\left(\scig\left(\cN_1\times\cN_2,(\nabla\Psi_1,\nabla\Psi_2)\right)\right)=\upsg\left(\scig\left(\cN_1,\nabla\Psi_1\right)\right) \star \upsg\left(\scig\left(\cN_2,\nabla\Psi_2\right)\right).
\]
\et

The following theorem is one of the most important fact in our reasoning. It allows us to simplify the distinguishing of the infinite-dimensional equivariant Conley indexes, significantly.

Let $\bH=\overline{\bigoplus_{n=0}^{\infty} \bH_n}$ be a representation of the compact Lie group G. Consider two functionals $\vp_1,\vp_2\in C^2_G(\bH,\bR)$ such that $\nabla\vp_i=Lx+\nabla K_i(x)$, where $\nabla K_i\in C^1_G(\bH,\bH)$ is completely continuous for $i=1,2$, which satisfy the conditions (B.1)-(B.3) described previously in subsection \ref{conleyIndexSection}.
Note that $T_x^{\perp} G(x)$ - a space orthogonal to the orbit - is a representation of the isotropy group $G_x$ and if $\vp$ is $G$-invariant then ${\vp}_{\mid T_{x}^{\perp} G(x)}$ is $G_x$-invariant.

\bt(\cite{PCHRYST2}, Theorem 2.4.3)\label{smash theorem for sci}
Let  $G(x_1), G(x_2)$ be isolated orbits of critical points of the potentials $\vp_1$ and $\vp_2$, respectively. Moreover, assume that $G_{x_1}=G_{x_2}(=H)$. If the pair $(G,H)$ is admissible and
$\upsh(\sci_H(\{x_1\},-\nabla \phi_1)) \neq \upsh(\sci_H(\{x_2\},-\nabla \phi_2)) \in U(H)$ where $\phi_i={\vp_i}_{\mid T_{x_i}^{\perp} G(x_i)}$ then 
\[\upsg(\sci_G(G(x_1),-\nabla \vp_1)) \neq  \upsg(\sci_G(G(x_2),-\nabla \vp_2)) \in U(G).
\]
\et

The proof of the theorem above is based on a concept of smash product over group. One can find more details in \cite{PCHRYST2}, especially Definition 2.4.2, Theorem 2.4.1 and Theorem 2.4.2.

\subsection{Equivariant splitting lemma}
\label{splitting}
\numsubsec 
 
Let $K$ be a compact Lie group and let $(\bV,\langle\cdot,\cdot\rangle)$ be an orthogonal Hilbert representation  of $K$  with an invariant scalar product $\langle\cdot,\cdot\rangle$. Moreover, assume that $\dim \bV^K < \infty$. Let $\Omega \subset \bV$ be an open and $K$-invariant neighborhood of $0\in\bV$.

 Consider a functional $\Psi\in C^2_K(\Omega,\bR)$of the form
\beq\label{psidef}
	\Psi(x)=\frac{1}{2}\langle Ax,x \rangle+\zeta(x),
\eeq
which satisfies the following assumptions
\be
\item [(F.1)] $A:\bV\to\bV$ is a $K$-equivariant self-adjoint linear Fredholm operator,
\item [(F.2)] $\ker A \subset \bV^K$,
\item [(F.3)] $\nabla\zeta:\bV\to \bV$ is a $K$-equivariant, compact operator,
\item [(F.4)] $\nabla\zeta(0)=0$ and $||\nabla^2\zeta(0)||= 0$,
\item [(F.5)] $0 \in \Omega$ is an isolated critical point of $\Psi$.
\ee

Note that the kernel $\ker A$ and the image $\im A$ are orthogonal representations of $K$. Moreover, $\ker A$ is  finite dimensional and trivial representation of $K.$ Since $A$ is self-adjoint, $\bV=\ker A \oplus \im A.$ Put $x=(u,v)$, where $u \in \ker A$ and $v \in \im A.$

The following theorem (called \emph{splitting lemma}) provides the existence of equivariant homotopy which allows us to study the product (splitted) flow $(\nabla\varphi(u),Av)$, where $u\in\ker A,\,v\in\im A$ instead of the general $\Psi(x)=\frac{1}{2}\langle Ax,x \rangle+\zeta(x)$. The proof of this theorem one can find in \cite{PCHRYST2} (Theorem 2.5.2).

\bt\label{splittinglemma} Suppose that a functional $\Psi\in C^2_K(\Omega,\bR)$ is defined by formula \eqref{psidef} and satisfies assumptions $(F.1)$--$(F.5)$. Then, there exists $\varepsilon_0>0$ and $K$-equivariant  gradient homotopy $\nabla \cH:(B_{\varepsilon_0}(\ker A) \times B_{\varepsilon_0}(\im A)) \times [0,1] \to \bV$ satisfying the following conditions:
\begin{enumerate}
\item $\nabla \cH((u,v),t)=Av-\nabla \xi_t(u,v)$, for $t\in [0,1]$, where $\nabla \xi_t=\nabla \xi(\cdot,t)$ and $\nabla \xi:\bV\times [0,1]\to\bV$ is  compact and  $K$--equivariant.
\item $(\nabla \cH)^{-1}(0)\cap (B_{\varepsilon_0}(\ker A)\times B_{\varepsilon_0}(\im A))\times[0,1]=\{0\}\times[0,1]$ i.e. $0$ is an isolated critical point of $\nabla\cH(\cdot,t)$ for any $t\in[0,1]$.
\item $\nabla \cH((u,v),0)=\nabla \Psi (u,v)$.
\item There exists an $K$--equivariant, gradient mapping $\nabla \varphi:B_{\varepsilon_0}(\ker A) \to \ker A$ such that $\nabla \cH((u,v),1)=(\nabla\varphi(u),Av),$ for all $(u,v) \in B_{\varepsilon_0}(\ker A) \times B_{\varepsilon_0}(\im A).$
\end{enumerate}
\et

\br\label{grad_varphi}
The homotopy $\cH$ is given by
\[
\cH((u,v),t)=\frac{1}{2}\langle Av,v\rangle +\frac{1}{2} t(2-t)\langle Aw(u),w(u)\rangle
		+ t\zeta (u,w(u))+(1-t)\zeta(u,v+tw(u)),
\]
Moreover, from the proof of Theorem \ref{splittinglemma} follows that the potential $\varphi: B_{\varepsilon_0}(\ker A) \to \bR$ is given by $\varphi(u)=\Psi(u,w(u))$, where $w:B_{\varepsilon}(\ker A)\to\im A\cap \bV^{K}$ is $K$-equivariant, see Remark 2.5.1 in \cite{PCHRYST2}.
\er

\br 
Note that we don't assume that $\ker A\neq\{0\}$. In the case of trivial kernel the homotopy given in Theorem \ref{splittinglemma} provides a linearization of functional.
\er

\section{Variational formulation for Hamiltonian systems}\label{variationa}
\numsec
Recall that we are interested in the existence of periodic solutions with any period of the system \eqref{hamsys}. In order to find them we are going to study $2\pi$-periodic solutions of the parameterized system
\beq\label{hamsys2}\tag{HS-P}
\dot z(t)=\lambda J \nabla H(z(t)),
\eeq 
which are in one-to-one correspondence to $2\pi\lambda$-periodic solutions of the system \eqref{hamsys}.

To prove the existence of solutions of the Hamiltonian system \eqref{hamsys2} we are going to the study critical points of a corresponding functional.

Define the Sobolev space of $2\pi$-periodic $\brtwon$-valued functions
\[
\bE:=\{z(t)=a_0+\sum_{k=1}^{\infty} a_k\cos(kt)+b_k\sin(kt)\,:\, a_i,b_i\in\brtwon,\, \sum_{k=1}^{\infty} k(|a_k|^2+|b_k|^2)<\infty\}.
\]
Then $\bE=\overline{\bE_0\oplus\bigoplus_{k=1}^{\infty} \bE_k}$ where $\bE_0=\brtwon$ is a subspace of constant functions and $\bE_k=\{a\cos(kt)+b\sin(kt)\,:\,a,b\in\brtwon\}$. Moreover  $\bE_k=\bE_k^-\oplus\bE_k^+$ for $k\geq 1$, where
\[
\bE_k^-=\{a\cos (kt)+Ja\sin(kt)\,:\,a\in\brtwon\},
\] 
\[
\bE_k^+=\{a\cos (kt)-Ja\sin(kt)\,:\,a\in\brtwon\}.
\]
The space $\bE$ with inner product given by
\beq\label{D:scalarProduct}
\langle z,\tilde{z}\rangle_{\bE}:=2\pi a_0\cdot \tilde{a}_0+\pi\sum_{k=1}^{\infty} k(a_k\cdot \tilde{a}_k+b_k\cdot \tilde{b}_k),
\eeq
where $\cdot$ denotes the standard scalar product, is a Hilbert space usually denoted by $\bH^{1/2}(\sone,\brtwon)$. Since we consider $\brtwon$ as a unitary representation of the compact Lie group $\Gamma$, $\bE$ is a unitary $G=\Gamma\times\sone$-representation with the action given by

\beq\label{actionDefinition}
G \times \bE \ni ((\gamma,e^{i \theta}),z(t)) \to \gamma z(t+\theta) 
\eeq
and $\bE_k$ is a unitary $G$-representation for any $k\geq 1$.
Indeed, 
\[
a_k\cos(k(t+\theta))+b_k\sin(k(t+\theta))=\begin{bmatrix}
\cos(k\theta)\cdot Id_{2N} &\sin(k\theta)\cdot Id_{2N} \\ -\sin(k\theta)\cdot Id_{2N}&\cos(k\theta)\cdot Id_{2N}
\end{bmatrix}
\begin{bmatrix}
a_k \cos(kt)\\ b_k\sin(kt)
\end{bmatrix}
\]
and therefore the action proposed in \eqref{actionDefinition} is given on $\bE_k$ by the product of unitary matrices $\begin{bmatrix}
\gamma&0 \\ 0 & \gamma
\end{bmatrix}$ and $\begin{bmatrix}
\cos(k\theta)\cdot Id_{2N}&\sin(k\theta)\cdot Id_{2N} \\ -\sin(k\theta)\cdot Id_{2N} & \cos(k\theta)\cdot Id_{2N}
\end{bmatrix}$.

\br 
Since we are going to study the Hamiltonian system $\eqref{hamsys2}$ is a neighborhood of the the orbit $G(z_0)$ of critical points, without loss of generality we can assume that Hamiltonian $H$ satisfies the following growth restriction
\beq\label{growthRestriction}
|\nabla H(z)|\leq a_1+a_2 |z|^s\qquad\text{for some }a_1,\,a_2>0,\quad s\in[1,\infty).
\eeq
Indeed, we may choose $\bar{H}$ such that $\nabla\bar{H}$ is bounded (i.e. $s=1$) and $\bar{H}(z)=H(z)$ in a neighborhood of the orbit $G(z_0)$.
\er

It is known (see \cite{MAWI}) that periodic solutions of the system \eqref{hamsys2} are in one to one correspondence with $\sone$-orbits of critical points of a  potential $\Phi : \bE \times (0,\infty) \to \bR$ of a class $C^1$ defined by
\beq\label{PhiDef0}
\Phi(z,\lambda)=\frac{1}{2} \langle Lz,z\rangle_{\bE}+K_{\lambda}(z),
\eeq
where
\beq\label{PhiDef}
\langle Lz,z\rangle_{\bE}=\int_0^{2\pi} J\dot{z}(t)\cdot z(t)\,dt,\qquad K_{\lambda}(z)=\int_0^{2\pi} \lambda H(z(t))\,dt.
\eeq
Note that $\Phi(\cdot,\lambda)$ acts on the subspace of constant functions $\bE_0$ as $\Phi_{\mid \bE_0\times\bR}(z,\lambda)=2\pi\lambda H(z)$.
Moreover, $L$ is given explicit on $z(t)=a_0+\sum_{k=1}^{\infty} a_k\cos(kt)+b_k\sin(kt)$ by
\beq\label{lineraFourierForm}
(Lz)(t)=\sum_{k=1}^{\infty} Jb_k \cos(kt)-Ja_k\sin(kt),
\eeq
see \cite{GOL2011}, the formula $(3.3)$. 
 
Since we consider $\brtwon$ as a unitary representation of a group $\Gamma$ and $H$ is $\Gamma$-invariant, the potential $\Phi$ is $\Gamma$-invariant. Moreover, it is $\sone$-invariant since it acts on $2\pi$-periodic functions.

Recall that since the Hamiltonian $H$ is $\Gamma$-invariant, the solutions of the system \eqref{hamsys2} form $\Gamma$-orbits i.e. if $z_0$ is a solution on \eqref{hamsys2} then $\gamma z_0$ solves \eqref{hamsys2} for any $\gamma\in\Gamma$. Therefore we are going to study $G=\Gamma\times\sone$-orbits of critical points of the corresponding $G$-invariant potential $\Phi$ i.e. we are interested in  solutions of the system
\beq\label{gradeq}
\nabla_z \Phi(z,\lambda)=0.
\eeq
Note that $\nabla_z \Phi(z,\lambda)=Lz+\nabla K_{\lambda}(z)$, $L$ is a linear, self-adjoint and $G$-equivariant operator and $\nabla K_{\lambda}(z)$ is completely continuous. Since $\ker L=\bE_0$ and $L_{\mid \bE_k^{\pm}}=\pm Id$, the conditions (B.1)--(B.3) given on the page \pageref{B-properties} are satisfied.

Let $z_0\in (\nabla H)^{-1}(0)$ and consider a linear Hamiltonian system
\beq\label{hamsyslin}\tag{HS-L}
\dot z(t)=\lambda J A(z(t)-z_0),
\eeq
which has a form of \eqref{hamsys2} with $H(z)=\frac 12 A(z-z_0)\cdot (z-z_0)$. The variational potential has the form
$
\Phi_L(z,\lambda)=\frac{1}{2} \langle Lz,z\rangle_{\bE}+\frac{1}{2} \langle B_{\lambda} (z-z_0),z-z_0\rangle_{\bE},
$
where
\beq\label{B-rieszForm}
\begin{split}
\langle &B_{\lambda} z,z'\rangle_{\bE}=\int_0^{2\pi} \lambda A z(t) \cdot z'(t) \, dt=\\
&\int_0^{2\pi} \left((\lambda A)a_0+\sum_{k=1}^{\infty} (\lambda A)a_k\cos (kt)+(\lambda A)b_k\sin(kt)\right)\cdot\left(a_0'+\sum_{k=1}^{\infty} a_k'\cos (kt)+b_k'\sin(kt)\right) \\
&= 2\pi (\lambda A)a_0\cdot a_0'+\pi\sum_{k=1}^{\infty} (\lambda A)a_k\cdot a_k'+(\lambda A)b_kb_k'.
\end{split}
\eeq
Taking into account the scalar product in $\bE$ given by \eqref{D:scalarProduct} and the formula \eqref{B-rieszForm} we obtain
\[
B_{\lambda} (z-z_0)=(\lambda A)(a_0-z_0)+\sum_{k=1}^{\infty} \frac{\lambda}k Aa_k \cos(kt)+\frac{\lambda}k Ab_k\sin(kt)
\]
and as a consequence
\[
\nabla_z\Phi_L(z,\lambda)=(L+B_{\lambda})(z-z_0)=\lambda A(a_0-z_0)+\sum_{k=1}^{\infty} \left(\frac{\lambda}k Aa_k +Jb_k\right)\cos(kt)+\left(\frac{\lambda}k Ab_k-Ja_k\right)\sin(kt).
\]
It means that $\nabla \Phi_L(z,\lambda)$ acts on $\bE_k=\{a\cos(kt)+b\sin(kt)\,:\,a,b\in\brtwon\}$ for $k\geq 1$ as a linear map
\beq\label{LinearDecomp}
T_{k,\lambda}(A)=\begin{bmatrix}
-\frac{\lambda}k A & -J \\ J & -\frac{\lambda}k A
\end{bmatrix}.
\eeq

\bl\label{linearLemma}
The linear equation \eqref{hamsyslin} possesses a non-constant $2\pi$-periodic solution if and only if $T_{k,\lambda}(A)$ is singular for some $k\geq 1$ and it holds true if $\lambda=
\frac{k}{\beta_r}$ where $i\beta_r\in\sigma(JA)$.
\el
\begin{proof}
Let $(z,\lambda)=(a_0+\sum_{k=1}^{\infty} a_k\cos (kt)+b_k\sin(kt),\lambda)\neq (const.,\lambda)$ be a critical point of $\Phi_L$ and let $k$ be such that $|a_k|^2+|b_k|^2\neq 0$. Then, in particular, $T_{k,\lambda}(A) (a_k,b_k)^T=0$ i.e. $T_{k,\lambda}(A)$ has a nontrivial kernel. 

It is easy to see that equation $T_{k,\lambda}(A)(a_k,b_k)^T=0$ has the form
\[
\left\{\begin{array}{rcl}
-\frac{\lambda}kAa_k & = & Jb_k \\
Ja_k & = & \frac{\lambda}k Ab_k \end{array}\right.
\]
which implies $JA(a_k-ib_k)=\frac{k}{\lambda}(b_k+ia_k)=\frac{ki}{\lambda}(a_k-ib_k)$ i.e. $\frac{ki}{\lambda}\in\sigma(JA)$.

\end{proof}

\section{Main Result}\label{S:main}
In this section we prove our main result of this paper i.e. the global bifurcation of periodic solutions of the system \eqref{hamsys} in the most general version. We emphasize our assumptions

\begin{enumerate}
\item[\textbf{(A1)}] $H:\brtwon\to\bR$ is a $\Gamma$-invariant Hamiltonian of the class $C^2$,
\item[\textbf{(A2)}] $z_0 \in \bR^{2N}$ is a~critical point of $H$ such that the isotropy group $\Gamma_{z_0}$ is trivial,
\item[\textbf{(A3)}] the orbit $\Gamma(z_0)$ is isolated in $(\nabla H)^{-1}(0)$,
\item[\textbf{(A4)}] $\pm i\beta_1,\ldots, \pm i\beta_m$, $0<\beta_m<\ldots<\beta_1$, $m\geq 1$ are the purely imaginary eigenvalues of $J\nabla^2 H(z_0)$,
\item[\textbf{(A5)}] $\deg(\nabla H_{\mid T^{\perp}_{z_0} \Gamma(z_0)},B(z_0,\epsilon),0) \neq 0$ for sufficiently small $\epsilon$,
\item[\textbf{(A6)}] $\beta_{j_0}$ is such that $\beta_j \slash \beta_{j_0}\not \in \bN$ for all $j\neq j_0$
\item[\textbf{(A7)}] $m^-\left(T_{1,\lambda}(\nabla^2 H(z_0))\right)$ changes at $\lambda=\frac 1{\beta_{j_0}}$ when $\lambda$ varies.
\end{enumerate}

\bt\label{T:general}
Under the assumptions \textbf{(A1)--(A7)} there exists a connected family of non-stationary periodic solutions of the system $\dot z(t)=J \nabla H(z(t))$ emanating from the stationary solution $z_0$ (i.e. with amplitude tending to 0) such that minimal periods of solutions in a small neighborhood of $z_0$ are close to $ 2\pi \slash \beta_{j_0}$.
\et

\br 
The assumption \textbf{(A7)} is very general and laborious to verify. We will change and simplify them in some specific cases. However, it does not follow directly from the structure of a Hamiltonian system in general situation as we obtained in a study of Newtonian systems, see \cite{PCHRYST}, the proof of Lemma 4.1.
\er

Let $z_0\in\brtwon$ be a critical point of the Hamiltonian $H$ such that the assumptions \textbf{(A1)--(A4)} are satisfied. From now we study variational reformulation \eqref{gradeq} of the parameterized Hamiltonian system \eqref{hamsys2}.  Then $z_0$ is a constant functions which solves the equation \eqref{gradeq} for any $\lambda\in (0,\infty)$ and the orbit $G(z_0)=\Gamma(z_0)$ consists of solutions of the equation \eqref{gradeq}. Therefore we put $\cT=G(z_0)\times (0,\infty)$ for the family of trivial solutions of the equation \eqref{gradeq} and $\cN=\{(z,\lambda)\in \bE \times (0,+\infty)\setminus \cT\,:\, \nabla_z\Phi(z,\lambda)=0\}$ is called a family of non-trivial solutions.

Denote by $\cC(z_0,\lambda_0)$ a connected component of the set $\overline{\cN}$ which contains the set $\{z_0\}\times\{\lambda_0\}$. 
\bdf\label{D:globBif}
We say that the orbit $G(z_0)\times\{\lambda_0\}$ is an orbit of \emph{global bifurcation} of solutions of the equation \eqref{gradeq} if the set $\cC(z_0,\lambda_0)$ is unbounded in $\bE\times (0,\infty)$ or $\left(\cC(z_0,\lambda_0)\cap \cT\right) \setminus \left(G(z_0)\times\lambda_0\right)\neq\emptyset$ i.e. $\cC(z_0,\lambda_0)$ coincide with the trivial family outside the orbit $G(z_0)\times\{\lambda_0\}$. 
\edf

The definition above does not depend on the choice of $z\in G(z_0)$. Indeed, if $z_1=g_1z_0$ then, taking into account an equivariancy of the equation \eqref{gradeq}, we obtain $\cC(z_1,\lambda_0)=g_1\cC(z_0,\lambda_0)$ i.e. the connected component of $\{z_1\}\times\{\lambda_0\}$  satisfies the same conditions as the connected component of $\{z_0\}\times\{\lambda_0\}$. In other words, global bifurcation from the orbit $G(z_0)\times\{\lambda_0\}$ provides the existence of solutions emanating from any point of the orbit. In fact, using the equivariant method we obtain the existence the bifurcation of the $G$-orbits of solutions. However, we are working with the bifurcation of single solutions (not orbits) to connect the main result of the paper to the original theorem of Lyapunov directly.

\br Note that the definition of global bifurcation implies that the set $\cC(z_0,\lambda_0)$ is not empty i.e. there is a family of solutions of the equation \eqref{gradeq} emanating from the orbit $G(z_0)\times\{\lambda_0\}$ at the point $\{z_0\}\times\{\lambda_0\}$. Therefore, to prove Theorem \ref{T:general} we have to show the existence of global bifurcation from the orbit $G(z_0)\times\{\lambda_0\}$ and to control the bifurcation level $\{\lambda_0\}$ to determine periods of bifurcating solutions. Finally, since the existence of bifurcation provides the convergence in the norm of Sobolev space $\bE=\bH^{1/2}(\sone,\brtwon)$, we have to prove that new periodic solutions tend to $\{z_0\}$ in the $L^{\infty}$-norm.
\er

Put $\Lambda=\{\frac{k}{\beta_j}\,:\, k\in\bN,\,i\beta_j\in\sigma(J\nabla^2_zH(z_0))\}$. In the theorem below we prove the necessary condition for the existence of bifurcation from the orbit $G(z_0)\times\{\lambda_0\}$.

\bt\label{necCond} (Necessary condition)
If $G(z_0)\times\{\lambda_0\}$ is an orbit of global bifurcation of solutions of the equation \eqref{gradeq} then $\ker \nabla^2_z\Phi(z_0,\lambda_0) \cap \overline{\bigoplus_{k=1}^{\infty} \bE_k}\neq\emptyset$ i.e. $\lambda_0 \in \Lambda.$
\et
\begin{proof}
By a reasoning given in the proof of Theorem 3.2.1 in \cite{PCHRYST2} we obtain $\ker \nabla^2_z\Phi(z_0,\lambda_0) \cap \overline{\bigoplus_{k=1}^{\infty} \bE_k}\neq\emptyset$. To complete the proof we have to prove that it implies $\lambda_0\in\Lambda$. The study of the kernel of $\nabla^2_z\Phi(z_0,\lambda_0)$ is equivalent to the study of the linearized system \eqref{hamsyslin} where $A=\nabla^2_zH(z_0)$. Therefore by the Lemma \ref{linearLemma} we obtain the thesis.      
\end{proof}

Choose $\lambda_0$ such that the necessary condition and assumptions \textbf{(A6), (A7)} are satisfied i.e. $\lambda_0=\frac{1}{\beta_{j_0}}\in\Lambda$ and put $\lpm=\frac{1\pm\varepsilon}{\beta_{j_0}}$ such that $\lambda_\pm>0$ and $[\lambda_-,\lambda_+]\cap\Lambda=\{\lambda_0\}.$ To prove the existence of global bifurcation we are going to apply the following theorem

\bt \label{sufficient} (Sufficient condition).
Under the assumptions above, if
\beq\label{suffCondition}
\upsg\left(\sci_{G}\left(G(z_0),-\nabla\Phi(\cdot,\lambda_+)\right)\right)\neq \upsg\left(\sci_{G}\left(G(z_0),-\nabla\Phi(\cdot,\lambda_-)\right)\right),
\eeq
then  $G(z_0)\times\{\lambda_0\}$ is an orbit of global bifurcation.
\et
\begin{proof}
The theorem above follows directly from the relation 
\[
\upsg(\scig(X,f)=\nabla_G-deg(f, \text{Int} X)
\]
(see \cite{BLGORY}, Theorem 3.10) and from a global bifurcation theorem for equivariant gradient degree (see \cite{GORY}, Theorem 3.3).
\end{proof}

Define $\bH\subset \bE$ by $\bH=T_{z_0}^{\perp} G(z_0)$. Recall that the space perpendicular to the orbit at $z_0$ is an $G_{z_0}$-representation. Since $z_0$ is a constant function and by the assumption \textbf{(A2)} $G_{z_0}=\{e\}\times\sone$, $\bH$ is an unitary $\sone$-representation. 

Put $\Psi_{\pm}:\bH\to\bR$ by $\Psi_{\pm}(z)=\Phi(z,\lpm)$. Note that since $\bH$ is an $\sone$-representation, the potenatial $\Psi_{\pm}$ is $\sone$-invariant. Moreover, $z_0$ is an isolated critical point of $\Psi_{\pm}$. Since $G(z_0)=\Gamma(z_0)\subset \bE_0$ we have the following decomposition
\[
\bH=T_{z_0}^{\perp} \Gamma(z_0)\oplus \overline{\bigoplus_{k=1}^{\infty} \bE_k}.
\]

In order to prove the main result of this paper we prove the existence of global bifurcation from the orbit $G(z_0)\times\{\lambda_0\}$ i.e. we need to prove formula \eqref{suffCondition}. In the theorem below we simplify this formula to the study of potentials defined on the orthogonal section $T_{z_0}^{\perp} G(z_0)$.

\bl\label{sectionConditionLemma}
Under the above assumptions if 
\beq\label{sectionCondition}
\upssone\left(\sci_{\sone}(\{z_0\},-\nabla \Psi_{+})\right)\neq \upssone\left(\sci_{\sone}(\{z_0\},-\nabla \Psi_{-})\right)
\eeq
then
\beq
\upsg\left(\sci_{G}\left(G(z_0),-\nabla\Phi(\cdot,\lambda_+)\right)\right)\neq \upsg\left(\sci_{G}\left(G(z_0),-\nabla\Phi(\cdot,\lambda_-)\right)\right).
\eeq
\el
\begin{proof}
Since the pair $(\Gamma\times\sone,\{e\}\times\sone)$ is admissible (because $\sone$ is abelian), see Definition \ref{admissible} and both $\nabla\Phi(\cdot,\lambda_-)$, $\nabla\Phi(\cdot,\lambda_+)$ are in the form of a compact perturbation of the same linear operator $L$, we can apply Theorem \ref{smash theorem for sci} to obtain the thesis directly.
\end{proof}

From now our goal is to prove formula \eqref{sectionCondition}. The next step is to transform a problem into the study of Conley indexes with simpler structure of flows.

We define $\widetilde{H}:T_{z_0}^{\perp} \Gamma(z_0) \to \bR$ by  $\widetilde{H}(z)=H(z+z_0)$ and $\widetilde{\Psi}_{\pm}:\bH\to\bR$ by $\widetilde{\Psi}_{\pm}(z)=\Psi_{\pm}(z+z_0)$. Since $\widetilde{\Psi}_{\pm\mid \bE_0}=2\pi\lpm\widetilde{H}$ and the orbits $G(z_0)\times\{\lambda_+\},\,G(z_0)\times\{\lambda_-\}$ do not satisfy the necessary condition for the existence of bifurcation we obtain $\ker\nabla^2\widetilde{\Psi}_{\pm|}(0)=\ker \nabla^2 \tilde{H}(0)$ so the kernel is independent on $\lambda_{\pm}$. Since $\nabla^2\tilde{\Psi}_{\pm}(0)$ is self-adjoint we are able to decompose
\[
\bH=\cN\oplus\cR=\ker\nabla^2\widetilde{\Psi}_{\pm}(0)\oplus \im\nabla^2\widetilde{\Psi}_{\pm}(0)
\]
independently on $\lambda$. We further decompose $\cR=\cR_0\oplus\cR_{\infty}$, where $\cR_0=\cR\cap\bE_0\subset \bH^{\sone}$ and $\cR_{\infty}=\overline{\bigoplus_{k=1}^{\infty} \bE_k}\subset\cR\subset\bH$. 
Note that the $\sone$-invariant  potential $\Pi_{\pm} : \cR_{\infty} \to \bR$ of the linear vector field $\nabla^2{\Psi_{\pm}}_{\mid \cR_{\infty}}(z_0)(z-z_0)$ is defined by $\ds \Pi_{\pm}(z)=\frac{1}{2} \langle \nabla^2{\Psi_{\pm}}_{\mid \cR_{\infty}}(z_0)(z-z_0),z-z_0\rangle_{\bE}.$

The next theorem simplifies the proof of formula \eqref{sectionCondition} to the study of Conley indexes of linear vector fields. In order to prove it we apply splitting lemma (Theorem \ref{splittinglemma}).

\bl\label{main-lemma0}
Under the above assumptions  
the formula \eqref{sectionCondition} holds true 
if and only if

 \beq\label{linearIndex}
 \upssone\left(\sci_{\sone}\left(\{z_0\},-\nabla\Pi_{+}\right)\right)\neq \upssone\left(\sci_{\sone}\left(\{z_0\},-\nabla\Pi_{-}\right)\right).
\eeq
\el 
\begin{proof}
It is clear that by the properties of Conley index we have 
\[
\sci_{\sone}(\{z_0\},-\nabla\Psi_{\pm})=\sci_{\sone}(\{0\},-\nabla\widetilde{\Psi}_{\pm}).
\]
Since we are going to apply splitting lemma (Theorem \ref{splittinglemma}), now we verify that $\widetilde{\Psi}_{\pm}$ satisfies conditions (F.1)--(F.5) given on the page \pageref{psidef} with $K=\sone$, $A=\nabla^2 \widetilde{\Psi}_{\pm}(0)$ and $\zeta_{\pm}(z)=\widetilde{\Psi}_{\pm}(z)-\langle\nabla^2 \widetilde{\Psi}_{\pm}(0)z,z\rangle_{\bE}$.
\begin{itemize}
\item [(F.1)] Since $\widetilde{\Psi}_{\pm}$ is $\sone$-invariant (it is the invariant $\Psi$ translated by $z_0\in\bE^{\sone}$) its hessian is $\sone$-equivariant. Moreover, a hessian is a self-adjont operator. By Theorem \ref{necCond} $\ker \nabla^2 \widetilde{\Psi}_{\pm}\subset \bE_0\cap \bH$ is finite dimensional, since $\lpm\notin\Lambda$.
\item [(F.2)] Similarly as above, $\ker \nabla^2 \widetilde{\Psi}_{\pm}\subset \bE_0\cap \bH=\bH^{\sone}$
\item [(F.3)] Since $\nabla \zeta_{\pm}(z)=\nabla\widetilde{\Psi}_{\pm}(z)-\nabla^2 \widetilde{\Psi}_{\pm}(0)z$ and both summands are compact and $\sone$-equivariant, $\nabla \zeta_{\pm}$ is also compact and $\sone$-equivariant.
\item [(F.4)] It is obvious due to formula given in (F.3).
\item [(F.5)] Since $\lpm\notin\Lambda$ i.e. the orbits $G(z_0)\times\{\lpm\}$ do not satisfy the necessary conditions for the existence of bifurcations, the orbit $G(z_0)$  is isolated in the set $(\nabla \Phi(\cdot,\lpm))^{-1}(0)$. Therefore $0\in\bH$ is an isolated critical point of $\tilde{\Psi}_{\pm}$.
\end{itemize}

\noindent Applying Theorem \ref{splittinglemma} (splitting lemma) and Theorem \ref{product formula} (product formula) we obtain
\[
\upssone\left(\sci_{\sone}(\{0\},-\nabla \widetilde{\Psi}_{\pm})\right)
=\upssone\left(\sci_{\sone}(\{0\},-\nabla\varphi_{\pm})\right)\star \upssone\left(\sci_{\sone}(\{0\},-\nabla^2 \widetilde{\Psi}_{\pm}(0)_{\mid \cR})\right),
\]
where $0=(0,0)\in \cN\oplus \cR $, $\varphi_{\pm} : B_{\varepsilon_0}(\cN) \to \bR$, $\varphi_{\pm}(u)=\widetilde{\Psi}_{\pm}(u,w(u))$ and $\nabla\varphi_{\pm}(u)$ is $\sone$-equivariant. 

 Since $\cR_0$ is an invariant space of the linear map $\nabla^2\widetilde{\Psi}_{\pm}(0)$ we are able to decompose the linear flow to obtain
\begin{multline*}
\upssone\left(\sci_{\sone}(\{(0,0)\},-\nabla^2 \widetilde{\Psi}_{\pm}(0)_{\mid \cR})\right)=
\\=\upssone\left(\sci_{\sone}(\{0\},-\nabla^2 \widetilde{\Psi}_{\pm}(0)_{\mid \cR_0})\right)\star \upssone\left(\sci_{\sone}(\{0\},-\nabla^2 \widetilde{\Psi}_{\pm}(0)_{\mid \cR_{\infty}})\right)
\end{multline*}
and combining the flows given on the $\bE_0\cap\bH$ we finally obtain
\begin{multline}\label{helpEq0}
\upssone\left(\sci_{\sone}(\{0\},-\nabla \widetilde{\Psi}_{\pm})\right)=\\
=\upssone\left(\sci_{\sone}(\{(0,0)\},(-\nabla\varphi_{\pm},-\nabla^2 \widetilde{\Psi}_{\pm}(0)_{\mid \cR_0}))\right)\star \upssone\left(\sci_{\sone}(\{0\},-\nabla^2 \widetilde{\Psi}_{\pm}(0)_{\mid \cR_{\infty}})\right).
\end{multline}
If we study the homotopy $\cH$ (see Theorem \ref{splittinglemma} and Remark \ref{grad_varphi}) acting on the subspace of constant function $\bE_0$ we obtain
\[
\begin{split}
\nabla\cH_{\mid \bE_0}((u,v),1)&=(-\nabla\varphi_{\pm}(u),-\nabla^2 \widetilde{\Psi}_{\pm}(0)_{\mid \cR_0}(v)),
\\
\nabla\cH_{\mid \bE_0}((u,v),0)&=\nabla\widetilde{\Psi}_{\pm\mid \bE_0}(u,v)=\widetilde{L}_{\mid\bE_0}(u,v)+\nabla\widetilde{K}_{\lpm\mid\bE_0}(u,v)=2\pi\lpm\nabla \widetilde{H}(u,v).
\end{split}
\]
By the homotopy invariance of the Conley index and since $\lambda_-,\,\lambda_+$ are both positive we have
\begin{multline*}
\upssone\left(\sci_{\sone}(\{(0,0)\},(-\nabla\varphi_{\pm},-\nabla^2 \widetilde{\Psi}_{\pm}(0)_{\mid \cR_0}))\right)=\\
=\upssone(\sci_{\sone}(\{0\},-\nabla \widetilde{H}))=\upssone(\sci_{\sone}(\{z_0\},-\nabla H_{\mid T_{z_0}^{\perp}\Gamma(z_0)})).
\end{multline*}
Note that the space $\bE_0$ is finite-dimensional and consists of constant functions (elements invariant on $\sone$ action), therefore
\begin{multline}\label{helpEq1}
\upssone(\sci_{\sone}(\{z_0\},-\nabla H_{\mid T_{z_0}^{\perp}\Gamma(z_0)}))=\chi_{\sone}(\ci_{\sone}(\{z_0\},-\nabla H_{\mid T_{z_0}^{\perp}\Gamma(z_0)}))=\\=\chi(\ci(\{z_0\},-\nabla H_{\mid T_{z_0}^{\perp}\Gamma(z_0)}))=\deg (\nabla H_{\mid T_{z_0}^{\perp}\Gamma(z_0)},B(z_0,\varepsilon))\cdot\bI\in U(\sone)
\end{multline}
for sufficiently small $\varepsilon>0$, where the last equality follows from Poincar\'e-Hopf theorem, see \cite{SRZEDNICKI}.
Since
\[
\upssone\left(\sci_{\sone}(\{0\},-\nabla^2 \widetilde{\Psi}_{\pm}(0)_{\mid \cR_{\infty}})\right)=\upssone\left(\sci_{\sone}(\{z_0\},-\nabla \Pi_{\pm})\right),
\]
we finally have
\beq\label{helpEq2}
\upssone\left(\sci_{\sone}(\{0\},-\nabla \widetilde{\Psi}_{\pm})\right)=\deg (\nabla H_{\mid T_{z_0}^{\perp}\Gamma(z_0)},B(z_0,\varepsilon))\cdot \upssone\left(\sci_{\sone}(\{z_0\},-\nabla \Pi_{\pm})\right)
\eeq
By the assumption \textbf{(A5)} $\deg(\nabla H_{\mid T^{\perp}_{z_0} \Gamma(z_0)},B(z_0,\epsilon),0) \neq 0$ and due to equation \eqref{helpEq2}
we obtain that the formula \eqref{sectionCondition} is equivalent to 
\[
\upssone\left(\sci_{\sone}\left(\{z_0\},-\nabla\Pi_{+}\right)\right)\neq \upssone\left(\sci_{\sone}\left(\{z_0\},-\nabla\Pi_{-}\right)\right)
\]
and the proof is completed.
\end{proof}

To verify formula \eqref{linearIndex} we are going to study equivariant Conley index and equivariant Euler characteristic by definitions. Note that the vector field $-\nabla\Pi_{\pm}:\bR_{\infty}\to\bR_{\infty}$ is linear and the decomposition $\cR_{\infty}=\overline{\bigoplus_{k=1}^{\infty} \bE_k}$ satisfies conditions (B.1)--(B.3) given on the page \pageref{B-properties}. Recall that $\lambda_0=\frac{1}{\beta_{j_0}}\in\Lambda$ where $i\beta_{j_0}\in\sigma(J\nabla^2 H(z_0))$ and $\lpm=\frac{1+\varepsilon}{\beta_{j_0}}$ is such that $[\lambda_-,\lambda_+]\cap\Lambda=\{\lambda_0\}$.

\br 
Note that the linearization of the variational functional $\Phi$ for the parameterized Hamiltonian system $\eqref{hamsys2}$ is equal to variational functional for the  linearized system $\eqref{hamsyslin}$ (we remove high order tenses in both cases). Therefore the action of the linear vector field $\nabla\Pi_{\pm}:\bR_{\infty}\to\bR_{\infty}$ is given on $\bE_k$ by 
\beq\label{LinearDecomp2}
T_{k,\lambda}(A)=\begin{bmatrix}
-\frac{\lambda}k A & -J \\ J & -\frac{\lambda}k A
\end{bmatrix}.
\eeq
where $A=\nabla^2 \widetilde{H}(0)=\nabla^2 H(z_0)$. For $k\to\infty$ we have $T_{k,\lambda}(A)\to \begin{bmatrix}
0& -J \\ J & 0
\end{bmatrix}$ i.e. $m^-(T_{k,\lambda})=2N$ for $k$ large enough, say for $k\geq k_0$.
\er

\bt\label{lastLemma}
Under the assumptions \textbf{(A1)}--\textbf{(A7)} of Theorem \ref{T:general}
\beq\label{linearDifference}
\upssone\left(\sci_{\sone}\left(\{z_0\},-\nabla\Pi_{+}\right)\right)\neq \upssone\left(\sci_{\sone}\left(\{z_0\},-\nabla\Pi_{-}\right)\right).
\eeq
\et
\begin{proof}
Since $T_{k,\lambda}$ is singular iff $\lambda=\frac{k}{\beta_j}$ (see Lemma \ref{linearLemma}), $[\lambda_-,\lambda_+]\cap\Lambda=\{\lambda_0\}=\{\frac{1}{\beta_{j_0}}\}$ and $\frac{1}{\beta_{j_0}}\neq\frac{k}{\beta_j}$ for any $k\in\bN$ and $\beta_j\neq\beta_{j_0}$ (see assumption \textbf{(A6)}), matrices $T_{k,\lambda}$ for $k\geq 2$ are nonsingular if $\lambda$ varies. 
Therefore the spectral decomposition of $\bE_k$ for $k\geq 2$ given by $-\nabla^2\Pi_{\pm}$ does not depend on $\lpm$ i.e
\[
\bE_k=\bE_{k,-}\oplus\bE_{k,+}
\]
but
\[
\bE_1=\bE_{1,\lpm,-}\oplus\bE_{1,\lpm,+}.
\]
As a consequence the spectra $\cE_-$, $\cE_+$ whose homotopy types are Conley indexes $\sci_{\sone}\left(\{z_0\},-\nabla\Pi_{-}\right),$ 
$\sci_{\sone}\left(\{z_0\},-\nabla\Pi_{+}\right)$ are of the same type $\xi=(\bE_{k,+})_{k=2}^{\infty}$. Define $\bP_n=\bigoplus_{k=2}^{n} \bE_{k,+}$.

Put $\cR_n=\bigoplus_{k=1}^{n} \bE_k$ and consider $\Pi_{\pm}^n=\Pi_{\pm\mid \cR_n}:\cR_n\to\bR$. By Remark \ref{stabilization of upsilon} we obtain
\beq
\upssone\left(\sci_{\sone}\left(\{z_0\},-\nabla\Pi_{\pm}\right)\right)=\left(\chi_{\sone}\left(S^{\bP_{n-1}}   \right)\right)^{-1}\star  \chi_{\sone}\left(\ci_{\sone}\left(\{z_0\},-\nabla\Pi^{n}_{\pm}\right)\right)
\eeq
for $n$ large enough. Now, to prove formula $\eqref{linearDifference}$ it is enough to show 
\[
\chi_{\sone}\left(\ci_{\sone}\left(\{z_0\},-\nabla\Pi^{n}_{-}\right)\right)\neq \chi_{\sone}\left(\ci_{\sone}\left(\{z_0\},-\nabla\Pi^{n}_{+}\right)\right).
\]
Since $-\nabla\Pi^n_{\pm}$ is a linear isomorphism Conley indexes are very simple i.e. 
\beq\label{helpFormula5}
\ci_{\sone}\left(\{z_0\},-\nabla\Pi^{n}_{\pm}\right)=S^{\bE_{1,\lpm,+}\oplus \bP_n}=S^{\bE_{1,\lpm,+}}\wedge S^{\bP_n}.
\eeq
By the assumption \textbf{(A7)} 
\[
\dim \bE_{1,\lambda_-,+}=m^-\left(T_{1,\lambda_-}(\nabla^2 H(z_0))\right)\neq m^-\left(T_{1,\lambda_+}(\nabla^2 H(z_0))\right)=\dim \bE_{1,\lambda_+,+}.
\]
Since $\bE_1$ is non-trivial $\sone$-representation by Remark \ref{propOfChi} we obtain
\beq\label{helpFormula6}
\chi_{\sone} \left(S^{\bE_{1,\lambda_-,+}}\right)\neq \chi_{\sone} \left(S^{\bE_{1,\lambda_+,+}}\right).
\eeq
Combining formulas \eqref{helpFormula5} and \eqref{helpFormula6} we finally obtain
\[
\begin{split}
\chi_{\sone}\left(\ci_{\sone}\left(\{z_0\},-\nabla\Pi^{n}_{-}\right)\right)&=\chi_{\sone}\left(S^{\bE_{1,\lambda_-,+}}\right)\star \chi_{\sone}\left(S^{\bP_n}\right)\neq \\
&\neq\chi_{\sone}\left(S^{\bE_{1,\lambda_+,+}}\right)\star \chi_{\sone}\left(S^{\bP_n}\right)=\chi_{\sone}\left(\ci_{\sone}\left(\{z_0\},-\nabla\Pi^{n}_{+}\right)\right)
\end{split}
\]
which completes the proof.
\end{proof}

\br\label{R:final1}
Theorem \ref{lastLemma} completes the proof of the existence of global bifurcation of solutions of the equation \eqref{gradeq} from the orbit $G(z_0)\times\{\lambda_0\}$. As a consequence we obtain the existence of connected branch of solutions of the system \eqref{hamsys} emanating from the stationary solutions $z_0$ with periods close to $2\pi\lambda_0=\frac{2\pi}{\beta_{j_0}}$. By the non-resonance condition for eigenvalues (i.e $\beta_j \slash \beta_{j_0}\not \in \bN$ for all $j\neq j_0$) we obtain $\frac{\lambda_0}{r}\notin\Lambda=\{\frac{k}{\beta_r}\,:\, k\in\bN,\,i\beta_r\in\sigma(J\nabla^2_zH(z_0))\}$ for any $r\in\bN$ and therefore there are no $\frac{2\pi\lambda_0}{r}$-periodic non-stationary solutions in a neighborhood of the orbit $\Gamma(z_0)$. Hence we can consider periods tending to $\frac{2\pi}{\beta_{j_0}}$ as minimal and the proof of Theorem \ref{T:general} is completed.
\er

\br\label{R:final2}
Applying bifurcation theory to the variational potential $\Phi:\bH^{1/2}(S^1,\brtwon)\times(0,\infty)\to\bR$ we prove the existence of family of critical points of $\Phi$ emanating from $\{z_0\}\times\{\lambda_0\}\in \bH^{1/2}(S^1,\brtwon)\times(0,\infty)$ in the norm of $\bH^{1/2}(S^1,\brtwon)$. Now we prove that corresponding periodic solutions of Hamiltonian system tend to $z_0\in\brtwon$ in $L^{\infty}$-norm. Let $z(t)=a_0+\sum_{k=1}^{\infty} a_k\cos(kt)+b_k\sin(kt)$ be a solution of \eqref{hamsys2} for $\lambda$ close to $\lambda_0$. Firstly,
\[
||z-z_0||^2_{L^2}=2\pi |a_0-z_0|^2+\pi\sum_{k=1}^{\infty} (|a_k|^2+|b_k|^2)\leq 2\pi |a_0-z_0|^2+\pi\sum_{k=1}^{\infty} k(|a_k|^2+|b_k|^2)=||z-z_0||^2_{\bH^1}.
\]
Under the condition \eqref{growthRestriction}, the map $z(t)\to\nabla H(z(t))$ is continuous from $L^2(\sone)$ to $L^2(\sone)$ (see Proposition B.1 in \cite{RAB5}). Let $\varepsilon>0$ and choose $0<\delta<\varepsilon$ such that $||z-z_0||_{L^2}\leq||z-z_0||_{\bH^1}< \delta$ implies $||\nabla H(z)||_{L^2}=||\nabla H(z)-\nabla H(z_0)||_{L^2}<\varepsilon$. Since $z$ is a solution of \eqref{hamsys2} we obtain
\[
||(z-z_0)'||_{L^2}=||\dot{z}||_{L^2}=||\lambda\nabla H(z)||_{L^2}\leq \lambda\varepsilon.
\]
Applying Sobolev inequality (see Proposition 1.1 in \cite{MAWI}) we obtain
\[
||z-z_0||^2_{L^{\infty}}\leq c||z-z_0||^2_{\bH^1}=c\left(||z-z_0||^2_{L^2}+||(z-z_0)'||^2_{L^2}\right)\leq c(1+\lambda^2)\varepsilon^2.
\]
Since $\lambda$ is bounded in the neighborhood of $\lambda_0$ the convergence of solutions $z$ to $z_0$ in the norm of $\bH^1$ implies the convergence in $L^{\infty}$ which completes the proof.

\er

\br 
The assumption \textbf{(A6)} was used only in the proof of Theorem \ref{lastLemma} i.e. in the last step of the proof of our main theorem. We are able to remove this assumption, but then in the proof of Theorem \ref{lastLemma} we need to study $\lambda$-depending decompositions of not only $\bE_1$ but any $\bE_k$ such that $\frac{1}{\beta_{j_0}}\neq\frac{k}{\beta_j}$ for some $j$. It will cause a complicated notation and the proof will be less readable. However, the change of $\dim \bE_{k,\lambda,+}$ when $\lambda$ varies we will obtain in the same way as for $k=1$.
Note that assumption \textbf{(A6)} is always satisfied for $j_0=1$ since $\beta_1$ is the maximum of $\beta_i$.

\er

\br Lets summarize the proof of Theorem \ref{T:general} in the steps.
\begin{enumerate}
\item By the change of variables we translate the equation \eqref{hamsys} into \eqref{hamsys2}.
\item We formulate the equation \eqref{hamsys2} as a variational problem \eqref{gradeq}.
\item We apply equivariant Conley index and equivariant Euler characteristic to provide the existence of global bifurcation of solutions of the equation \eqref{gradeq} from the orbit $G(z_0)\times\{\lambda_0\}$. From now we are going to prove formula \eqref{suffCondition} i.e. the change of the equivariant gradient degree at the level $\lambda_0$ satisfying the necessary condition.
\item To study the change of equivariant Conley index of the orbit we apply the method of orthogonal section, reducing the problem to formula \eqref{sectionCondition}.
\item Applying equivariant splitting lemma and the assumption \textbf{(A5)} we reduce formula \eqref{sectionCondition} to the linear case i.e. to formula \eqref{linearIndex}.
\item Finally we prove formula \eqref{linearIndex} computing equivariant Conley index by the definition.
\end{enumerate}
\er

\section{Corollaries}
In this section we study in which way is it possible to modify assumption \textbf{(A7)}. Moreover, we show that the results of this paper are generalizations of some versions of Lyapunov center theorem.

The following theorem was proven by Szulkin (\cite{SZULKIN}, Proposition 3.6)
\bt
Suppose that $A$ is symmetric and $i\beta_j$, $\beta_j>0$, is an eigenvalue of $JA$. Let $E_j$ be the eigenspace of $JA$ in $\bC^{2N}$ corresponding to $i\beta_j$ and $Z_j$ the invariant subspace of $JA$ in $\bR^{2N}$ corresponding to $\pm i\beta_j$. Then $m^-(T_{1,\lambda}(A))$ changes at $\lambda=1/\beta_j$ if and only if the following two equivalent conditions are satisfied:
\begin{enumerate}
\item $m^-(A_{\mid Z_j})\neq m^+(A_{\mid Z_j})$,
\item $m^-(-iJ_{\mid E_j})\neq m^+(-iJ_{\mid E_j})$.
\end{enumerate}
\et 
\noindent Due to the theorem above we are able to formulate new versions of assumption \textbf{(A7)}:
\begin{itemize}
\item[\textbf{(A7.1)}] $m^-(\nabla^2 H(z_0)_{\mid Z_{j_0}})\neq m^+(\nabla^2 H(z_0)_{\mid Z_{j_0}})$, 
\item[\textbf{(A7.2)}] $m^-(-iJ_{\mid E_{j_0}})\neq m^+(-iJ_{\mid E_{j_0}})$, 
\end{itemize}
where $E_j$ be the eigenspace of $J\nabla^2 H(z_0)$ in $\bC^{2N}$ corresponding to $i\beta_j$ and $Z_j$ the invariant subspace of $J\nabla^2 H(z_0)$ in $\bR^{2N}$ corresponding to $\pm i\beta_j$.

Note that if $\nabla^2 H(z_0)_{\mid Z_{j_0}}$ is a definite matrix then the condition \textbf{(A7.1)} is satisfied.  Therefore we put
\begin{itemize}
\item[\textbf{(A7.3)}] $\nabla^2 H(z_0)_{\mid Z_{j_0}}$  is definite.
\end{itemize}

\bt\label{T:cor1}
Under the assumptions \textbf{(A1)--(A6)} and one of the conditions \textbf{(A7.1)--(A7.3)} there exists a connected family of non-stationary periodic solutions of the system $\dot z(t)=J \nabla H(z(t))$ emanating from the stationary solution $z_0$ such that minimal periods of solutions in the small neighborhood of $z_0$ are close to $ 2\pi \slash \beta_{j_0}$.
\et

If we are not interested in the minimal period of new solutions but only in the study of its existence, the assumptions can be modified. The computation of invariant subspaces $Z_j$ we can change to the study of general invariant subspace of $J\nabla^2 H(z_0)$ associated to all the eigenvalues of the form $\pm i\beta_k$. Denoting this subspace by $Z$ we formulate a new condition.  

\begin{itemize}
\item[\textbf{(A7.4)}] $\nabla^2 H(z_0)_{\mid Z}$  is definite.
\end{itemize}

Under this condition the assumption \textbf{(A7.3)} is satisfied for some eigenvalue of $J\nabla^2 H(z_0)$ and we do not know it precisely. Therefore we exclude assumption \textbf{(A6)}. In the theorem below we prove the existence of periodic solutions of the system \eqref{hamsys} without information about their minimal periods. Under the reasoning above, it is clear that Theorem \ref{T:cor2} is a direct consequence of Theorem \ref{T:general}.
\bt\label{T:cor2}
Under the assumptions \textbf{(A1)--(A5)} and \textbf{(A7.4)} there exists a connected family of non-stationary periodic solutions of the system $\dot z(t)=J \nabla H(z(t))$ emanating from the stationary solution $z_0$ such that periods (not necessarily minimal) of solutions in the small neighborhood of $z_0$ are close to $ 2\pi \slash \beta_{j}$ where $i\beta_j$, $\beta_j>0$, is some eigenvalue of $J\nabla^2 H(z_0)$.
\et

Looking on the $T_{1,\lambda}(\nabla^2 H(z_0))$ from the other point of view, we study the limit behaviour of the Morse index $m^-\left(T_{1,\lambda}(\nabla^2 H(z_0))\right)$ as $\lambda\to 0$ and as $\lambda\to\infty$.

When $\lambda\to 0$, we have
\[
T_{1,\lambda}(\nabla^2 H(z_0)) \to \begin{bmatrix}
0 & -J \\ J & 0
\end{bmatrix}.
\]
The limit matrix is nonsingular (its eigenvalues are $\pm 1$, each of multiplicity $2N$), and therefore the Morse index $m^-\left(T_{1,\lambda}(\nabla^2 H(z_0))\right)$ is constant and equal to $2N$ on some interval $\lambda\in(0,\lambda_r)$.

On the other hand,
\[
\frac{1}{\lambda}T_{1,\lambda}(\nabla^2 H(z_0)) \to \begin{bmatrix}
-\nabla^2 H(z_0) & 0 \\ 0 & -\nabla^2 H(z_0)
\end{bmatrix}\quad\text{ for }\lambda\to \infty.
\]
This limit matrix is, in general, singular, since $\ker \nabla^2 H(z_0)$ can be nontrivial (see Remark \ref{eqgrad}). However, by Lemma \ref{linearLemma}, the matrix $T_{1,\lambda}(\nabla^2 H(z_0))$ is nonsingular except at a finite number of levels $\lambda$. As a consequence, all zero eigenvalues of the limit matrix appear as the limits of positive or negative eigenvalues of $\frac{1}{\lambda}T_{1,\lambda}(\nabla^2 H(z_0))$. Therefore, since $m^-\left(T_{1,\lambda}(\nabla^2 H(z_0))\right)=m^-\left(\frac{1}{\lambda}T_{1,\lambda}(\nabla^2 H(z_0))\right)$, for sufficiently large $\lambda$ we can bound the Morse index as follows:
\[
\begin{split}
2m^+(\nabla^2 H(z_0))&=2m^-(-\nabla^2 H(z_0))\leq
m^-\left(T_{1,\lambda}\left(\nabla^2 H(z_0)\right)\right)\\
&\leq 2\left(2N-m^+(-\nabla^2 H(z_0))\right)
=2\left(2N-m^-(\nabla^2 H(z_0))\right).
\end{split}
\]
Consequently, the condition $m^-\left(T_{1,\lambda}\left(\nabla^2 H(z_0)\right)\right) \neq 2N$ for large $\lambda$, which implies that the Morse index $m^-(T_{1,\lambda}(\nabla^2 H(z_0)))$ changes at some $\lambda\in (0,\infty)$, is guaranteed if $m^+(\nabla^2 H(z_0))> N$ or $m^-(\nabla^2 H(z_0)) > N$. Recall that the levels $\lambda$ where it can change is $\Lambda$ (see Lemma \ref{linearLemma}). Therefore the change of $m^-(T_{1,\lambda}(\nabla^2 H(z_0)))$ implies the existence of purely imaginary eigenvalue of $J\nabla^2 H(z_0)$.
As a consequence we can propose new condition
\begin{itemize}
\item[\textbf{(A7.5)}] $m^+(\nabla^2 H(z_0))> N$ or $m^-(\nabla^2 H(z_0)) > N$
\end{itemize}
and we are able to formulate the next theorem without assumption \textbf{(A4)}. 

\bt\label{T:cor3}
Under the assumptions \textbf{(A1),(A2),(A3),(A5)} and \textbf{(A7.5)} there exists a connected family of non-stationary periodic solutions of the system $\dot z(t)=J \nabla H(z(t))$ emanating from the stationary solution $z_0$ such that periods (not necessarily minimal) of solutions in the small neighborhood of $z_0$ are close to $ 2\pi \slash \beta_{j}$, where $i\beta_j$, $\beta_j>0$, is some eigenvalue of $J\nabla^2 H(z_0)$.
\et

Below we present in which way the theorems presented above generalize classical Lyapunov center theorem and an analogous theorem for Hamiltonian systems that has been proved by Dancer and Rybicki \cite{DARY}. Moreover, two symmetric version of the Lyapunov center theorem proposed in \cite{PCHRYST} and \cite{PCHRYST2} are generalized in this paper.

\bt\label{T:DARY} (\cite{DARY}, Theorem 3.3. (reformulated))
Consider an equation $\dot z(t)=J \nabla H(z(t))$, where $H:\brtwon\to\bR$ is of the class $C^2$. Let $z_0\in\brtwon$ be an isolated critical point of $H$. Let $i\beta_0$ $\beta_0>0$ be an eigenvalue of $J\nabla^2 H(z_0)$. If $\deg(\nabla H, B(z_0,\epsilon),0) \neq 0$ for sufficiently small $\epsilon$ and $m^-\left(T_{1,\lambda}(\nabla^2 H(z_0))\right)$ changes at $\lambda=\frac 1{\beta_{0}}$ when $\lambda$ varies, then there exists a connected family of periodic solution of the Hamiltonian system emanating from $z_0$.
\et
\begin{proof}
This theorem follows directly from Theorem \ref{T:general} if we consider trivial group $\Gamma=\{e\}$. In this case $T_{z_0} \Gamma(z_0)=\brtwon+z_0$.
\end{proof}
Consider a Newtonian (second-order) system
\beq\label{newsys2}\tag{NS}
\ddot{q}(t)=-\nabla U(q(t)),
\eeq
where $U:\bR^N\to\bR$ is $\Gamma$-invariant potential of the class $C^2$, $\Gamma$ acts orthogonally on $\bR^N$ and $q_0\in(\nabla U)^{-1}(0)$. If we substitute $r=\dot{q}$ we can reformulate the second-order system \eqref{newsys2} to the first-order system
\beq\label{newtonToHamilton}
\left\{
\ba{l}
\dot{q}(t)=r(t),\\
\dot{r}(t)=-\nabla U(q(t)),
\ea
\right.
\eeq
which can be considered as a Hamiltonian system with $H:\brtwon
\to\bR$ defined by $H(z)=H(q,r)=\frac 12 r^2+U(q)$. An action of $\Gamma$ on $\brtwon$ induced by action on $\bR^N$ is diagonal i.e $\gamma(q,r)\to(\gamma q,\gamma r)$. It is easy to verify that this action is symplectic, so $\Gamma$ acts unitary on $\brtwon$. Moreover, $z_0=(q_0,r_0)=(q_0,\dot{q}_0)=(q_0,0)$ (since we consider $q_0$ as a constant function) is a critical point of $H$.
We see that $J\nabla^2 H(z_0)=\left[\ba{rr} 0 & \nabla^2 U(q_0) \\ -I & 0 \ea \right]$. The easy block--form of the matrix $J\nabla^2 H(z_0)$ lets us to observe a bijective correspondence between positive eigenvalues of $\nabla^2 U(q_0)$ and the pairs of purely imaginary eigenvalues of $J\nabla^2 H(z_0)$. In fact, if $\beta^2\in\sigma(\nabla^2 U(q_0))$ then $\pm i\beta \in\sigma(J\nabla^2 H(z_0))$. 
Taking into account the above reasoning, the following theorems are consequences of Theorem \ref{T:general}.

\bt[Symmetric Lyapunov center theorem, \cite{PCHRYST}]\label{SLCT}

Let $U:\Omega\to\bR$ be a $\Gamma$-invariant potential of the class $C^2$ and $q_0 \in \Omega$. Assume that
\begin{enumerate}
\item  $q_0$ is a critical points of the potential $U$,
\item  $\dim \ker \nabla^2 U(q_0) = \dim \Gamma(q_0),$
\item  the isotropy group $\Gamma_{q_0}$ is trivial,
\item  $\sigma(\nabla^2 U(q_0)) \cap (0,+\infty) = \{\beta_1^2,\ldots, \beta_m^2\}$ and $m\geq 1$.
\end{enumerate}
Then for any $\beta_{j_0}$ such that $\beta_j \slash \beta_{j_0}\not \in \mathbb{N}$ for all $j\neq j_0$ there exists a sequence $(q_k(t))$ of periodic solutions of the system \eqref{newsys2}  with minimal period tending to $2 \pi \slash \beta_{j_0}$ such that in any open neighborhood of the orbit $\Gamma(q_0)$ there is an element of the sequence $(q_k(t))$.
\et 

\bt[Symmetric Lyapunov center theorem for minimal orbit, \cite{PCHRYST2}]\label{SLCTM}
Let $U:\Omega\to\bR$ be a $\Gamma$-invariant potential of the class $C^2$ and $q_0 \in \Omega$. Assume that
\begin{enumerate}
\item $q_0$ is a minimum of potential $U$
\item the orbit  $\Gamma(q_0)$ is isolated in $(\nabla U)^{-1}(0)$,
\item the isotropy group $\Gamma_{q_0}$ is trivial,
\item $\sigma(\nabla^2 U(q_0)) \cap (0,+\infty) = \{\beta_1^2,\ldots, \beta_m^2\}$, $\beta_1>\beta_2>\ldots>\beta_m>0$ and $m\geq 1$.
\end{enumerate}
Then for any $\beta_{j_0}$ such that $\beta_j \slash \beta_{j_0}\not \in \bN$ for $j\neq j_0$ there exists a sequence $(q_k(t))$ of periodic solutions of the system \eqref{newsys2} with a sequence of minimal periods $(T_k)$ such that $\mathrm{dist} (\Gamma(q_0),q_k([0,T_k]))\to 0$ and $T_k\to 2 \pi \slash \beta_{j_0}$ as $k\to\infty$.
\et 
\begin{proof}
Note, that the assumptions \textbf{(A1)}--\textbf{(A4)} are satisfied directly due to statements of the theorems above.

Firstly, we check that the assumption \textbf{(A7)} is always satisfied for Newtonian systems. Since $\nabla^2 H(z_0)=\begin{bmatrix}
\nabla^2 U(q_0) & 0 \\ 0 & I
\end{bmatrix}$ and the matrix $\nabla^2 U(q_0)$ is orthogonally diagonalizable (say by $D\in O(N)$) then the symplectic matrix $\bar{D}=\begin{bmatrix} D & 0 \\ 0 & D\end{bmatrix}$ diagonalize the hessian $\nabla^2 H(z_0)$ and we are able to simplify the form of $T_{1,\lambda}(\nabla^2 H(z_0))$ as follows
\[
\begin{bmatrix}
\bar{D} & 0 \\ 0 & \bar{D} 
\end{bmatrix}\cdot
\begin{bmatrix}
-\lambda \nabla^2 H(z_0) & -J \\ J & -\lambda \nabla^2 H(z_0)
\end{bmatrix}\cdot \begin{bmatrix}
\bar{D}^T & 0 \\ 0 & \bar{D}^T
\end{bmatrix}=\begin{bmatrix}
-\lambda \bar{D}\nabla^2 H(z_0)\bar{D}^T & -J \\ J& -\lambda \bar{D}\nabla^2 H(z_0)\bar{D}^T
\end{bmatrix}
\]
where $\bar{D}\nabla^2 H(z_0)\bar{D}^T=\diag (\eta_1,\eta_2,\ldots,\eta_N,1,\ldots,1)$ and $\eta_1,\eta_2,\ldots,\eta_N$ are the eigenvalues of $\nabla^2 U(q_0)$ (not necessarily different). 
Further, we apply the permutation of the basis $(1,2,\ldots,4N)\to(1,N+1,2N+1,3N+1,2,N+2,2N+2,3N+2,\ldots,N,2N,3N,4N)$ to transform our matrix to $\diag(A_1,A_2,\ldots,A_N)$, where $A_i=\begin{bmatrix}-\lambda\eta_i & 0 & 0 & 1 \\ 0 & -\lambda & -1 & 0 \\ 0 & -1 & -\lambda\eta_i & 0 \\ 1 & 0 & 0 & -\lambda\end{bmatrix}$. The characteristic polynomial of $A_i$ has the form $W_{A_i}(t)=\left((\lambda\eta_i+t)(\lambda+t)-1\right)^2$ and has $4$ negative roots for $\lambda^2=\lambda_+^2=\frac{1+\varepsilon}{\eta_i}$ and $2$ negative roots $\lambda^2=\lambda_-^2=\frac{1-\varepsilon}{\eta_i}$, where $\eta_i$ is positive. Therefore $m^-(T_{1,\lambda_+}(\nabla^2 H(z_0))=m^-(T_{1,\lambda_-}(\nabla^2 H(z_0)))+2mult(\eta_i)$ for sufficiently small $\varepsilon$, so the assumption \textbf{(A6)} is satisfied automatically for any positive eigenvalue $\beta_i^2$ of $\nabla^2 U(q_0)$ when we study Newtonian system translated into the Hamiltonian one.

To complete the proofs of theorems we have to verify assumption \textbf{(A5)} in both cases.
\begin{itemize}\label{tempTemp}
\item In Theorem \ref{SLCT} we assume that the orbit is non-degenerate i.e. $\dim \ker \nabla^2 U(q_0) = \dim \Gamma(q_0)$ and $T_{q_0} \Gamma(q_0)=\ker \nabla^2 U(q_0)$. Therefore $z_0$ is non-degenerate critical point of $\nabla H_{\mid T^{\perp}_{z_0} \Gamma(z_0)}$ i.e. $\nabla^2 H_{\mid T^{\perp}_{z_0} \Gamma(z_0)}$ is an isomorphism. In such case $\deg(\nabla H_{\mid T^{\perp}_{z_0} \Gamma(z_0)},B(z_0,\epsilon),0)=\pm 1$.
\item In Theorem \ref{SLCT} we assume that the orbit $\Gamma(q_0)$ consists of minima of $U$ and is isolated in critical points of $U$. Therefore $z_0=(q_0,0)$ is an isolated minimum of $H_{\mid T^{\perp}_{z_0} \Gamma(z_0)}$. However, it is known that Brouwer degree of minimum equals $1$, see \cite{RAB4}.
\end{itemize}
Theorems \ref{SLCT} and \ref{SLCTM} are given with the original thesis but in fact they are directly related to the thesis of Theorem \ref{T:general}, see remarks below Definition \ref{D:globBif}.

\end{proof} 

\section{An application}
In this section we apply our abstract results to the study of the quasi-periodic motions of a satellite near the geostationary orbit of an oblate spheroid with rotational symmetry. Note that the Earth is flattened and therefore the study of gravitation potential of such bodies has crucial role in the design of missions of satellites.

A gravitational potential of an oblate spheroid has a general form
\[
U_G(r,\theta)=-G\frac{E}{r}\left(1-\sum_{n=2}^{\infty}\left(\frac Rr\right)^n J_n P_n(\cos\theta)\right),
\]
where $r$ is a distance from the center of mass of the spheroid, $\theta$ - deviation from the axis of rotation, $G$ - gravitational constant, $E$ is the mass of the spheroid and $R$ is its equatorial radius, $(J_n)$ is the sequence of coefficients realted to the spherical harmonics and $P_n$ denotes the n-th Legendre polynomial, see \cite{LOWRIE} for the details.
In the case of axial symmetry the dominating term is $J_2$, so called dynamical form-factor, which is directly related to the body's flattening and for the oblate body $J_2$ is positive. For the Earth $J_2=1.0826359\cdot 10^{-3}$. 

We are going to study the motions under approximate potential
\[
U(r,\theta)=-\frac{GE}{r}\left(1-\frac{J_2 R^2}{r^2}P_2(\cos\theta)\right).
\]
By the choose of the units we may assume $GE=1$. Moreover, $P_2(x)=\frac{1}{2}\left(3x^2-1\right)$ and by the change of coordinates to the axially symmetric cylindrical ones we obtain
\beq\label{examplePotential}
V(r,z)=-\frac{1}{d}\left(1-\frac{c}{d^2}\left(3\frac {z^2}{d^2}-1\right)\right)=-\frac{1}{d}-\frac{c}{d^3}+\frac{3cz^2}{d^5},
\eeq
where $c=\frac 12 R^2J_2>0$, $d=\sqrt{r^2+z^2}$.

Assume that axially symmetric and oblate planet is rotating with an angular velocity $\omega$. We study the move of the satellite in the gravity field of this planet without influence of other bodies. Denote by $q_1,q_2,q_3$ coordinates of the satellite in a frame rotating with an angular velocity $\omega$ (the frame fixed with planet), where the axis $q_3$ is the axis of rotation and symmetry of the planet and denote by $p_1,p_2,p_3$ the corresponding momenta. The equation of motion is generated by the Hamiltonian $H$ of the form:
\beq\label{E:hamApp}
H(q_1,q_2,q_3,p_1,p_2,p_3)=\frac 12 (p_1^2+p_2^2+p_3^2)+\omega(q_1p_2-q_2p_1)+V(r,q_3),
\eeq
where $r=\sqrt{q_1^2+q_2^2}$ and $V$ is given in \eqref{examplePotential}, see \cite{H&F}. Note that this Hamiltonian is $\sone$-invariant where the symplectic action is given by
\[
\sone\times\bR^6 \ni (e^{i\theta},(q_1,q_2,q_3,p_1,p_2,p_3))\to (\begin{pmatrix}
\cos \theta & -\sin\theta&0 \\ \sin\theta& \cos\theta &0 \\ 0 &0&1
\end{pmatrix}\begin{pmatrix}
q_1\\q_2\\q_3
\end{pmatrix}, \begin{pmatrix}
\cos \theta & -\sin\theta&0 \\ \sin\theta& \cos\theta &0 \\ 0 &0&1
\end{pmatrix}\begin{pmatrix}
p_1\\p_2\\p_3
\end{pmatrix}).
\]
Non-zero equilibria of the Hamiltonian system $\dot{z}(t)=J\nabla H(z(t))$ describe a motion of a satellite along geostationary orbits. We apply Theorem \ref{T:cor3} to prove the existence of periodic solutions in a nearby of any equilibrium. Since the coordinates frame is rotating, we obtain the quasi-periodic motions of the satellite in a neighborhood of the geostationary orbit. We are interested in geostationary circular orbit so we assume $r>0$ 

Firstly, we have to find critical points of $H$. 
\[
\left\{\ba{l}
 H'_{q_1} (q_1,q_2,q_3,p_1,p_2,p_3)=\omega p_2 +V'_{r}(r,q_3)\frac{q_1}{r}, \\
 H'_{q_2} (q_1,q_2,q_3,p_1,p_2,p_3)=-\omega p_1 +V'_{r}(r,q_3)\frac{q_2}{r}, \\
  H'_{q_3} (q_1,q_2,q_3,p_1,p_2,p_3)=V'_{z}(r,q_3), \\
  H'_{p_1} (q_1,q_2,q_3,p_1,p_2,p_3)=p_1-\omega q_2, \\
  H'_{p_2} (q_1,q_2,q_3,p_1,p_2,p_3)= p_2+\omega q_1, \\ 
  H'_{p_3} (q_1,q_2,q_3,p_1,p_2,p_3) = p_3.
 \ea\right.
\]
Therefore, critical points of $H$ need to satisfy 
\[
\left\{\ba{l}
\left(r\omega^2-V'_r (r,q_3)\right)q_1=0, \\
\left(r\omega^2-V'_r (r,q_3)\right)q_2=0, \\
V'_z(r,q_3) = 0, \\
p_1=\omega q_2, \\
p_2= - \omega q_1, \\
p_3 = 0.

\ea\right.
\]
Since $r=\sqrt{q_1^2+q_2^2}>0$ by the first two equations we have
\beq\label{appEq1}
r\omega^2=V'_r (r,q_3)=\frac{r}{d^3}+\frac{3cr}{d^5}-\frac{15cq_3^2r}{d^7} \Rightarrow \omega^2 d^7=d^4+3cd^2-15cq_3^2.
\eeq
Further, by the third equation
$0=V'_z(r,q_3)=\frac{q_3}{d^3}+\frac{9cq_3}{d^5}-\frac{15cq_3^3}{d^7}=\frac{q_3}{d^7}\left(d^4+9cd^2-15cq_3^2\right)=\frac{q_3}{d^7}\left(w^2d^7+6cd^2\right)$. Since $\omega,c,d>0$ we obtain $q_3=0$. As a consequence the equation \eqref{appEq1} has a form
\[
\omega^2d^5-d^2-3c=0.
\]
Since $c>0$, by Descartes rule of signs there exists exactly one positive root of this equation, say $d_0$. It means that there exist one $SO(2)$ orbit of critical points of $H$ i.e. $SO(2)(Q)$ where $Q=(d_0,0,0,0,-\omega d_0,0)$. The point $Q$ from this orbit is chosen such that $r=\sqrt{q_1^2+q_2^2}=q_1=d_0$. This orbit is obviously isolated in $(\nabla H)^{-1}(0)$. To apply theorem \ref{T:cor3} we compute the Hessian $\nabla^2 H(Q)$. We have
\[
\nabla^2 H(Q)=\begin{bmatrix}
V''_{r,r}(d_0,0) & 0 & V''_{r,z}(d_0,0) & 0 &\omega &0 \\
0 & \omega^2 &0&-\omega&0&0\\
V''_{r,z}(d_0,0) &0 &V''_{z,z}(d_0,0) &0&0&0\\
0&-\omega&0&1&0&0 \\
\omega&0&0&0&1&0\\
0&0&0&0&0&1
\end{bmatrix}
\]
and
\beq\label{appEq2}
\ba{l}
V''_{r,r}(d_0,0)=-\frac{2}{d_0^3}-\frac{12c}{d_0^5}<0, \\
V''_{z,z}(d_0,0)=\frac{1}{d_0^3}+\frac{9c}{d_0^5}>0, \\
V''_{r,z}(d_0,0)=0.
\ea
\eeq
The Hessian is obviously degenerate (see Remark \ref{eqgrad}). One can see that it possesses eigenvalues $1+\omega^2$ (with eigenvector $[0,1,0,-1/\omega,0,0]^T]$) and $1$ (with eigenvector $[0,0,0,0,0,1]^T$). Denote be $\lambda_1,\lambda_2,\lambda_3$ the other three eigenvalues. If we compute the characteristic polynomial $w(t)$ of $\nabla^2 H(Q)$ its coefficient of the term $t$ (which is the additive inverse of the product of eigenvalues different from the one zero-eigenvalue we have already know) equals
\[
-(1+\omega^2)\left(-\omega^2 V''_{z,z}(d_0,0)+V''_{r,r}(d_0,0)V''_{z,z}(d_0,0)-(V''_{r,z}(d_0,0))^2\right)
\]
and substituting formulas \eqref{appEq2} we obtain
\[
\lambda_1\lambda_2\lambda_3= \left(\frac{1}{d_0^3}+\frac{9c}{d_0^5}\right)\left(-\omega^2 -\frac{2}{d_0^3}-\frac{12c}{d_0^5}\right)<0.
\]
Hence one or three of $\lambda_i$ are negative. Moreover, since $V''_{r,z}(d_0,0)=0$, the vector $[0,0,1,0,0,0]^T$ is an eigenvector of $\nabla^2 H(Q)$ corresponding to the positive eigenvalue $V''_{z,z}(d_0,0)$, so we may assume $\lambda_3=V''_{z,z}(d_0,0)>0$. Therefore exactly one of $\lambda_i$ is negative and the Hessian $\nabla^2 H(Q)$ has exactly four positive eigenvalues, i.e. $m^+(\nabla^2 H(Q))=4>3=N$. It means that the assumption \textbf{(A7.5)} is satisfied. Moreover, the kernel of this Hessian is one-dimensional which provides that the orbit $SO(2)(Q)$ is non-degenerate. Therefore the assumption \textbf{(A5)} is also satisfied (see the reasoning in the last paragraph of the previous section on the page \pageref{tempTemp}).
To summarize, all assumptions of Theorem \ref{T:cor3} are satisfied. It provides the existence of periodic solutions in a nearby of an equilibrium $Q$ in the rotating frame. These solutions correspond to a motion in a neighborhood of the geostationary orbit. 


\subsection*{Acknowledgements} 
I would like to thanks prof. S{\l}awomir Rybicki for the fruitful discussions on the topic of this article and prof. Andrzej Maciejewski for the proposition of physical-motivated example.

The author was partially supported by the National Science Centre, Poland (Grant No. 2017/25/N/ST1/00498).

\bibliographystyle{abbrv}
\bibliography{20261106bibliography-arXiv}

 \end{document}